\documentclass[12pt]{article} 
 
\usepackage{amssymb} 
\usepackage{a4} 
\newcommand{\sect}[1]{\setcounter{equation}{0}\section{#1}}

\newcommand{\app}[1]{\setcounter{section}{0} 
\setcounter{equation}{0} \renewcommand{\thesection} 
{\Alph{section}}\section{#1}} 
 
\newcommand{\A}{{\cal A}}

\def\1{{\bf 1}} 
\def\b#1{{\mathbb #1}} 
\def\c#1{{\cal #1}} 
\def\R{{\cal R}\,} 
 
\newcommand{\tl}{\,\triangleleft} 

\def\cocross{{>\!\!\!\triangleleft\,}} 
 
\def\id{\mbox{id\,}} 
\def \uqg{\mbox{$U_q{\/\mbox{\bf g}}$ }}

\def\nn{\nonumber \\} 
\newcommand{\be}{\begin{equation}} 
\newcommand{\ee}{\end{equation}} 
\newcommand{\bea}{\begin{eqnarray}} 
\newcommand{\eea}{\end{eqnarray}} 
\newcommand{\ba}{\begin{array}} 
\newcommand{\ea}{\end{array}} 
\newtheorem{prop}{Proposition}

\def\sq{\mbox{\rlap{$\sqcap$}$\sqcup$}} 
\newenvironment{proof}[1]{\vspace{5pt}\noindent{\bf Proof #1}\hspace{6pt}}%
{\hfill\sq} 
\newcommand{\bp}{\begin{proof}} 
\newcommand{\ep}{\end{proof}\par\vspace{10pt}\noindent} 
\begin{document} 
\title{Quantum group covariant (anti)symmetrizers,
$\varepsilon$-tensors, vielbein, Hodge map and Laplacian}

\author{        Gaetano Fiore  
\footnote{Work partially supported by the European Commission RTN Programme
HPRN-CT-2000-00131 and by MIUR} \\\\ 
         \and 
        Dip. di Matematica e Applicazioni, Fac.  di Ingegneria\\  
        Universit\`a di Napoli, V. Claudio 21, 80125 Napoli 
\\ 
         and\\
        I.N.F.N., Sezione di Napoli,\\ 
        Complesso MSA, V. Cintia, 80126 Napoli 
        } 
\date{}

\maketitle
\abstract{$GL_q(N)$- and $SO_q(N)$-covariant deformations of the
completely symmetric/antisymmetric projectors with an arbitrary
number of indices are explicitly constructed as polynomials
in the braid matrices. The precise relation between the 
completely antisymmetric projectors and the completely 
antisymmetric tensor is determined. Adopting the $GL_q(N)$- and
$SO_q(N)$-covariant differential calculi on the corresponding quantum 
group covariant noncommutative spaces $\b{C}_q^N$, $\b{R}_q^N$,
we introduce a generalized notion of vielbein basis (or ``frame''),
based on differential-operator-valued 1-forms. We then 
give a thorough definition of a  $SO_q(N)$-covariant
$\b{R}_q^N$-bilinear Hodge map acting 
on the bimodule of differential forms on $\b{R}_q^N$, introduce 
the exterior coderivative and show that the Laplacian
acts on differential forms exactly as in the undeformed case,
namely it acts on each component as it does on functions.}  

\vfill 
\noindent 
Preprint 04-5 Dip. Matematica e Applicazioni, Universit\`a di Napoli\\ 
DSF/ 08-2004 
\newpage

\section{Introduction and preliminaries}

The noncommutative geometry program \cite{Con,DubKerMad90} and the related program
of generalizing the concept of symmetries through quantum
groups \cite{Dri85-86,Wor87-89,FadResTak89} and quantum group covariant 
noncommutative spaces (shortly: quantum spaces) \cite{Man,FadResTak89} has
found a widespread interest in the mathematical and theoretical physics
community over the past two decades for its potential applications both in
fundamental and applied
 phsyics.   In order to make either program
powerful on a specific model it is important to reproduce as many of the tools
available in the corresponding undeformed (commutative) geometry model (if
any) as possible. 
The scope of the present work is to revisit and/or solve a number of
related technical issues, left (partially ot totally) untreated or unsolved
in the literature, regarding the quantum groups $H=U_qgl(N),U_qso(N)$
of the classical series  \cite{Dri85-86,FadResTak89}, the noncommutative
spaces $\b{C}_q^N,\b{R}_q^N$ \cite{Man,FadResTak89} on which they act,
and the quantum group covariant differential calculi
\cite{PusWor89,WesZum90,CarSchWat91}  on the latter.

As known, the braid matrix $\hat{R}$ of $H$ \cite{FadResTak89}
is a $N^2 \times N^2$ matrix, $H$-covariant deformation of
the permutation matrix $P$. $H$-covariant (anti)symmetrizers
${\cal P}^{\pm}$ of 2-tensors arise from the projector 
decomposition of $\hat R$, or equivalently can be expressed as
(first or second degree) polynomials in $\hat R$.
In analogy with the undeformed case, $H$-covariant (anti)symmetrizers  ${\cal
P}^{\pm,l}$  of  $l$-tensors ($l\ge 2$) are expected (see e.g.
\cite{Fio93,fiothesis})  to be polynomials
in $\hat{R}_{12},...,\hat{R}_{(l\!-\!1)l}$, the matrices obtained as tensor
products of $\hat{R}$ with $l\!-\!2$ copies of the $N\times N$ unit matrix. In
section \ref{comple} we find a very compact and manageable recursive relation,
through which these polynomials are determined. For $H=U_qso(N)$ this is in
agreement with the much longer recursive relation found in Ref.
\cite{HecSch99}\footnote{We thank the authors of  \cite{HecSch99} for calling
our attention to their paper, which we didn't know, after
the appearance of the first version of the present work on the
electronic arXive.}. In section \ref{eps} we
 recall or prove properties
of the $H$-covariant $\varepsilon$-tensor
\cite{Lyu87,Gur88,FadResTak89,Fio94,Maj95} and
 determine precisely
 its
relation with the antisymmetric projectors ${\cal P}^{-,l}$.
 In section
\ref{vielb} we introduce an open-minded  generalization of the
 notion
\cite{DimMad96}  of vielbein (or ``frame'') basis of
 1-forms on
$\b{C}_q^N,\b{R}_q^N$; we modify the approach adopted
 for $\b{R}_q^N$  in
\cite{CerFioMad00}, in that we allow the matrix 
 trasforming the basis
$\{dx^i\}$ into   
 the vielbein to have as entries 
 differential operators,
rather than functions\footnote{As  
 a by-product some unpleasent 
aspects of the vielbein of \cite{CerFioMad00} 
disappear. Incidentally, this change of 
attitude should allow to introduce a frame basis also  
for other quantum spaces, notably $q$-Minkowski.}. 
In section \ref{hodge} we introduce a thorough and  
consistent definition 
of a bilinear, $U_qso(N)$-covariant 
Hodge map, exterior coderivative and Laplacian
acting on differential forms on  $\b{R}_q^N$.

The projector decomposition of the
$H$-covariant braid matrix $\hat{R}$ reads 
\be
\ba{ll}
\hat R = q{\cal P}^+ - q^{-1}{\cal P}^-,\qquad\qquad\qquad\quad &\mbox{if }
H=U_qgl(N),      \qquad     \\[12pt]
\hat R = q{\cal P}^+ - q^{-1}{\cal P}^- + q^{1-N}{\cal P}^t,       
  \qquad\quad &\mbox{if } H=U_qso(N).     \qquad     
\ea                                        \label{projectorR} 
\ee
${\cal P}^-$
is the corresponding deformation of the antisymmetric 
projector. In (\ref{projectorR})$_1$ the matrix ${\cal P}^+$ is the
$U_qgl(N)$-covariant 
deformation of the symmetric projector, in (\ref{projectorR})
it is the $U_qso(N)$-covariant  deformation of the symmetric trace-free
projector, while ${\cal P}^t$ is the trace projector.
Thus they satisfy the equations 
\be 
{\cal P}^{\alpha}{\cal P}^{\beta} = {\cal P}^\alpha \delta^{\alpha\beta}, \qquad 
\sum_{\alpha}{\cal P}^{\alpha} = \1_{N^2},            \label{projector1} 
\ee 
where $(\1_{N^2})^{ij}_{hk}=\delta^i_h\delta^j_k$, 
$\alpha,\beta = -,+$ in the $H=U_qgl(N)$ case and $\alpha,\beta = -,+,t$ in the
$H=U_qso(N)$ case.  $\hat R$ is a symmetric matrix, and therefore also the
projectors are:
\be
\hat R^T=\hat R,\qquad\qquad {\cal P}^{\alpha}{}^T={\cal
P}^{\alpha}.\label{sym}  
\ee
The braid matrix fulfills the equation   
\be 
f(\hat R_{12})\,\hat R_{23}\,\hat R_{12}  
=\hat R_{23}\,\hat R_{12}\,f(\hat R_{23})                     \label{braid2} 
\ee 
for any rational function $f(t)$ in one variable such that 
the spectrum of $f(\hat R)$ has no poles, in particular 
for $f(\hat R)=\hat R,\hat R^{-1},{\cal P}^{\alpha}$.  
Here we have used the conventional matrix-tensor notation 
$\hat R_{12} = \hat R \otimes \1_N$,  
$\hat R_{23} =  \1_N\otimes \hat R $, 
where $\1_N$ denotes the $N\times N$ unit matrix.

In the $H=U_qso(N)$ case the ${\cal P}^t$ projects on a
one-dimensional sub-space and   can be written in the form  
\be 
{\cal P}^t{}_{kl}^{ij} = \frac 1{Q_N} g^{ij}g_{kl}           \label{Pt} 
\ee 
where the $N \times N$ matrix $g_{ij}$  
is a $U_qso(N)$-isotropic tensor,  
deformation of the ordinary Euclidean metric, which
will be given in (\ref{defgij}), and \cite{CarSchWat91}
$$
Q_N\equiv g^{lm}g_{lm}\!=\!\frac{(1\!+\!q^{2\!-\!N})(q^N\!-\!1)}{q^2-1}
\!=\!
\left(q^{1\!-\!\frac N2}\!+\!q^{\frac N2\!-\!1}\right)\left[\frac N2\right]_q.
$$
Here and in the sequel  we use the ``$q$-deformed numbers''
\be
[y]_q\!:=\!\frac {q^y\!-\!q^{-y}}{q\!-\!q^{-1}},\qquad
\qquad  y_{q^{\pm 2}}\!:=\!\frac {q^{\pm 2y}\!-\!1}{q{\pm 2}\!-\!1}=
q^{\pm (l\!-\!1)}[y]_q.
                     \label{defnum}
\ee
The metric and the braid matrix
satisfy the   relations~\cite{FadResTak89} 
\be 
g_{il}\,\hat R^{\pm 1}{}^{lh}_{jk} =  
\hat R^{\mp 1}{}^{hl}_{ij}\,g_{lk}, \qquad 
g^{il}\,\hat R^{\pm 1}{}_{lh}^{jk} =  
\hat R^{\mp 1}{}_{hl}^{ij}\,g^{lk}.                          \label{gRrel} 
\ee 
Indices will be lowered and raised using $g_{ij}$ and its inverse $g^{ij}$, 
e.g. 
$$ 
\partial^i:=g^{ij}\partial_j \qquad \qquad 
x_i:=g_{ij}x^j. 
$$ 

By taking powers of either decomposition (\ref{projectorR}) 
one can can express the projectors as polynomials in $\hat R$. One finds
\bea
{\cal P}^{\pm}&=&\frac{q^{\mp 1}\1_{N^2}\pm \hat R}{q+q^{-1}}\qquad\qquad
\qquad\qquad\qquad\mbox{if } H=U_qgl(N),\qquad       \label{projectorpm1}\\
{\cal P}^{\pm}&=&\frac{q^{\mp 1}\1_{N^2}\pm \hat R
-(q^{\mp 1}\pm q^{1\!-\!N}){\cal P}^t}{q+q^{-1}}\qquad
\qquad\mbox{if } H=U_qso(N).\qquad                 \label{projectorpm2}
\eea

\medskip 
  
The deformed algebras 
of functions on the two quantum  spaces  are called ``algebra
of functions on the quantum  hyperplane $\b{C}_q^N$'' and ``algebra
of functions on the quantum  Euclidean space $\b{R}_q^N$'' respectively
($h:=\ln q$ plays the role of deformation parameter);
we shall denote either one by $F$.   $F$ is essentially the unital
associative  algebra over $\b{C}[[h]]$ generated by $N$ elements $x^i$ (the
cartesian  ``coordinates'') modulo the relations  
(\ref{xxrel}) given below. 
The corresponding $H$-covariant differential calculi
\cite{WesZum90,CarSchWat91} are defined introducing  
the invariant exterior derivative $d$, satisfying nilpotency 
and the Leibniz rule $d(fg)=dfg+fdg$, and imposing the covariant 
commutation relations 
 (\ref{xxirel}) between the $x^i$ and the differentials 
$\xi^i:=dx^i$.  
 Partial derivatives are introduced through the  
decomposition 
$d=:\xi^i\partial_i$. All the other commutation relations 
are derived by consistency. The complete list is 
given by 
\bea 
&& {\cal P}^-{}^{ij}_{hk}x^hx^k=0, \label{xxrel}\\ 
&& x^h\xi^i=q\hat R^{hi}_{jk}\xi^jx^k,\label{xxirel}\\ 
&& (\1_{N^2}-{\cal P}^-)^{ij}_{hk}\xi^h\xi^k=0,\label{xixirel}\\ 
&& {\cal P}^-{}^{ij}_{hk}\partial_j\partial_i=0, \label{ddrel}\\ 
&& \partial_i x^j = \delta^j_i+q\hat R^{jh}_{ik} x^k\partial_h,  
                                                 \label{dxrel}\\ 
&& \partial_h\xi^i=q^{-1}\hat R^{-1}{}^{ik}_{hj}\xi^j
\partial_k.                                     \label{dxirel}  
\eea
We shall call ${\cal DC}^*$ (differential calculus algebra 
on $\b{R}_q^N$) the unital associative  
algebra  over $\b{C}[[h]]$ generated by  
$x^i,\xi^i,\partial_i$ modulo these relations. 
We shall denote by  
$\bigwedge^*$ (exterior algebra, or algebra of exterior forms) 
the graded unital subalgebra generated by the $\xi^i$ alone, with 
grading $\natural \equiv$the degree in $\xi^i$, and by $\bigwedge^p$  
(vector 
space of exterior $p$-forms) the component with grading $\natural =p$,  
$p=0,1,2...$. 
Each $\bigwedge^p$ carries an irreducible representation 
of $H$, and its dimension is the binomial 
coefficient $N\choose{p}$ \cite{FadResTak89,Fio94}, exactly as in the $q=1$ 
(i.e. undeformed) case; in particular there are no forms with 
$p>N$, and $\mbox{dim}(\bigwedge^N)={N\choose{N}}=1$, therefore 
$\bigwedge^N$ carries the singlet representation of $H$.
 
We shall endow ${\cal DC}^*$ with the same grading $\natural $, and call 
${\cal DC}^p$ its component with grading $\natural =p$. The elements of 
${\cal DC}^p$ can be considered differential-operator-valued $p$-forms. 
 
 We shall denote by $\Omega^*$ (algebra of differential forms) 
the graded unital subalgebra generated by the $\xi^i,x^i$, with 
grading $\natural $, and by $\Omega^p$  
(space of differential $p$-forms) its component with grading $p$;  
by definition $\Omega^0=F$ itself. 
Clearly both $\Omega^*$ and $\Omega^p$ are 
$F$-bimodules. 
 
 We shall denote by ${\cal H}$ (Heisenberg algebra) 
the unital subalgebra generated by the $x^i, \partial_i$. Note that 
by definition ${\cal DC}^0={\cal H}$, and that both 
${\cal DC}^*$ and ${\cal DC}^p$ are ${\cal H}$-bimodules. 
Finally, we shall call $F'$ the unital associative
algebra generated by the $\partial_i$ alone.
  
In the $H=U_qso(N)$-covariant case the elements  
$$ 
r^2:= x \cdot x=g_{kl} x^k x^l, \qquad \qquad  
\Box:= \partial\cdot\partial= 
g^{kl} \partial_l \partial_k=g_{kl}\partial^k \partial^l 
$$  
are $U_qso(N)$-invariant and respectively generate the centers 
of $F,F'$. 

The $H$-covariance of the differential calculus implies
that ${\cal DC}^*$ is a (right, in our conventions) 
$H$-module algebra.
All the information on the algebras ${\cal DC}^*,H$  
and the right action of the Hopf algebra $H$ on ${\cal DC}^*$
can be encoded in the  
cross-product algebra ${\cal DC}^*\cocross H$. 
We recall that this is $H\otimes {\cal DC}^*$ 
as a vector space, 
and so we denote as usual $g\otimes a$ simply by $ga$; 
that $H\1_{{\cal DC}^*}$, $\1_H{\cal DC}^*$ are subalgebras 
isomorphic to $H, {\cal DC}^*$, and so we omit to  
write either unit $\1_{{\cal DC}^*},\1_H$ 
whenever multiplied by non-unit elements;   
that for any $a\in{\cal DC}^*$, $g\in H$ the product fulfills 
\be 
a g=g_{(1)}\, (a\tl g_{(2)}).                \label{crossprod} 
\ee 
Here $\Delta(g)=g_{(1)}\otimes  g_{(2)}$ denotes the coproduct
of $g$ in Sweedler notation.
${\cal DC}^* \cocross H$ is a $H$-module algebra  itself, if we 
extend $\tl$ on $H$ as the adjoint action, namely as 
$$ 
h\tl g= Sg_{(1)}\, h\, g_{(2)}. 
$$ 
In view of (\ref{crossprod}), 
this formula will correctly reproduce the action also on 
the elements of ${\cal DC}^*$, and therefore on 
{\it any} element $a\in {\cal DC}^*\cocross H$. 
The elements $\sigma^i$, with $\sigma^i=x^i,\xi^i,\partial^i$, 
transform with the $N$-dimensional representation $\rho$
of $U_qsl(N)$ or $U_qso(N)$ respectively:
\be
\sigma^i\tl \, g=\rho^i_j(g)\sigma^j        \label{fundrep}
\ee

The above scheme applies also to the Hopf algebra
$H=\widetilde{U_qso(N)}$, which is the central extension of
$H=U_qso(N)$ obtained by adding
a central and primitive element $\eta$ generating dilatations 
of elements of ${\cal H}$,
\be
x^i\eta=(\eta+1)x^i, \qquad\quad 
\xi^i\eta=(\eta+1)\xi^i, \qquad\quad 
\partial^i\eta=(\eta-1)\partial^i.
\ee
We shall call $\eta$ also the generator of dilatations of
$U_qgl(N)$.

\medskip
One can introduce an alternative $H$-covariant differential calculus
replacing $q\hat R$ by $(q\hat R)^{-1}$ in the defining relations
(\ref{xxirel}-\ref{dxirel}).
The corresponding objects $\hat\xi^i,\hat\partial_i$ can be 
realized as suitable ``functions'' of  $x^j,\xi^j,\partial_j$ 
\cite{Ogi92}.

\sect{Completely (anti)symmetric projectors}
\label{comple}

The projectors 
${\cal P}^{\pm,l}\equiv \Vert{\cal P}^{\pm,l}{}^{i_1\ldots i_l}_{j_1\ldots
j_l}\Vert$ project the tensor product of  $l$ copies of the
$N$-dimensional representation of $H$ to the $l$-fold completely
symmetric/antisymmetric irreducible representation $V^{\pm,l}_N$ of $H$ 
therein contained. They
are uniquely characterized by the following properties
\bea 
&&{\cal P}^{\pm,l}{\cal P}^{\pi}_{m(m\!+\!1)}=
\delta^{\pm}_{\pi}{\cal P}^{\pm,l}=  
{\cal P}^{\pi}_{m(m\!+\!1)}{\cal P}^{\pm,l},  \label{Plproj1}  \\[8pt]
&&\left({\cal P}^{\pm,l}\right)^2={\cal P}^{\pm,l},   \label{Plproj2}  \\[8pt] 
&&\mbox{tr}_{1\ldots l}\!\left({\cal
P}^{\pm,l}\right) =\mbox{dim}(V^{\pm,l}_N),   \label{Plproj3}     
\eea
where $\pi=-,+$ in the $U_qsl(N)$- and $\pi=-,+,t$
in the $U_qso(N)$-covariant case respectively, $m=1,...,l\!-\!1$ and by
${\cal P}^{\pi}_{m(m\!+\!1)}$ we have denoted the matrix  
acting as ${\cal P}^{\pi}$ on the $m$-th, $(m\!+\!1)$-th indices and  
as the identity on the remaining ones. Eq. (\ref{Plproj3}) 
guarantees that ${\cal P}^{\pm,l}$ act as the identity (and not
as proper projectors) on $V^{\pm,l}_N$; for both $H=U_qgl(N),U_qso(N)$
$\mbox{dim}\!(V^{-,l}_N)\!=\!{
N\choose l}$, whereas for $H=U_qgl(N)$
$\mbox{dim}(V^{+,l}_N)={N\!-\!1\!+\!l\choose N\!-\!1 }$, 
for $H=U_qso(N)$
$\mbox{dim}(V^{+,l}_N)={N\!-\!1\!+\!l\choose N\!-\!1 }-{N\!-\!3\!+\!l\choose
N\!-\!1 }$.  In the appendix we prove

\begin{prop} The projectors ${\cal P}^{\pm,l\!+\!1}$ can be expressed
as polynomials in $\hat R_{12},...,\hat R_{(l\!-\!1)l}$ 
through either recursive relation
\bea
{\cal P}^{\pm,l\!+\!1}&=&
{\cal P}^{\pm,l}_{12...l}M^{\pm,l\!+\!1}_{l(l\!+\!1)}{\cal P}^{\pm,l}_{12...l},
\label{ansatz1} \\
&=& {\cal P}^{\pm,l}_{2...(l\!+\!1)}
M^{\pm,l\!+\!1}_{12}{\cal P}^{\pm,l}_{2...(l\!+\!1)}, \label{ansatz2} 
\eea
where ${\cal P}^{\pm,l}_{1...l}\equiv {\cal P}^{\pm,l}\otimes \1_N $,
${\cal P}^{\pm,l}_{2...(l\!+\!1)}\equiv \1_N\otimes {\cal P}^{\pm,l}$,
$M^{\pm,l\!+\!1}_{l(l\!+\!1)}= \1_N^{\otimes^{l\!-\!1}}\otimes
M^{\pm,l\!+\!1}$ etc., and   
\be
\ba{l}
M^{\pm,l\!+\!1} \!=\!\frac 1{[l\!+\!1]_q}
\!\left[q^{\mp l}\1_{N^2}\!\pm\![l]_q\hat R\right]\qquad\qquad\qquad\quad 
\mbox{if }H=U_qsl(N),\\[8pt] M^{\pm,l\!+\!1}  \!=\! \frac 1{[l\!+\!1]_q}
\!\left[q^{\mp l}\1_{N^2}\!\pm\![l]_q\hat R\!+\!
\frac{Q_N(q^{\pm 2}\!-\!1)[l]_q}{q^{\pm
1}\!\mp\!q^{N\!-\!1\!\pm\!2l}} 
{\cal P}^t\right]\quad \mbox{if }H=U_qso(N). 
\ea
                                  \label{Ml}
\ee
As a consequence, they are symmetric, $({\cal P}^{\pm,l})^T={\cal
P}^{\pm,l}$, and if $H=U_qso(N)$
\be
{\cal P}^{\pm,l}{}^{i_1...i_l}_{j_1...j_l} g^{j_1k_1}...g^{j_lk_l}
=g_{i_1j_1}...g_{i_lj_l}{\cal P}^{\pm,l}{}^{j_l...j_1}_{k_l...k_1}, 
\ee
and the same with the matrix $g^{ij}$ replaced by
its transpose $g^{ji}$.
\label{symmetrizers}
\end{prop}

In Ref. \cite{Fio93} we explicitly determined as examples just 
${\cal P}^{\pm,3}$  for $H=U_qso(N)$. In Ref. \cite{HecSch99}
longer recursive relations for ${\cal P}^{\pm,l}$ with arbitrary $l$
in the case $H=U_qso(N)$ were given; the Ansatz adopted there
was of the type 
${\cal P}^{\pm,l\!+\!1}=B^{\pm,l\!+\!1}{\cal P}^{\pm,l}_{12...l}$. The
unknown $N^{l\!+\!1}\times N^{l\!+\!1}$  matrices $B^{\pm,l\!+\!1}$
were explicitly determined
to be rather long polynomials in $\hat R_{12},...,\hat R_{(l\!-\!1)l}$. To go
from our formula to  theirs one just needs to set
$B^{\pm,l\!+\!1}={\cal P}^{\pm,l}_{12...l}M^{\pm,l\!+\!1}_{l(l\!+\!1)}$;
to go from their formula to ours one has to multiply both sides by
${\cal P}^{\pm,l}_{12...l}$ from the left and do a straightforward
calculation using (\ref{Plproj1}), (\ref{projectorR})$_2$.

{\bf Remark.} One can easily check that in the $H=U_qgl(N)$ case the deformed 
(anti)symmetric projectors ${\cal P}^{\pm,l}$
can be obtained from the polynomials 
giving the undeformed (anti)symmetric projectors in terms of the permutators
$P_{m(m\!+\!1)}$ by replacing the latter respectively with $\pm q^{\pm 1}\hat
R_{m(m\!+\!1)}$, and readjusting the  normalizations.

\sect{Properties of the $H$-covariant $\varepsilon$-tensors} 
\label{eps}

In our convention indices  $i,j,...$ take the following values:
\bea
&&i=1,2,\ldots,N \qquad\qquad\qquad\qquad\mbox{if }H=U_q gl(N), \\
&&i=-n,\ldots,-1,0,1,\ldots,n\qquad\quad\mbox{if }H=U_qso(2n\!+\!1), \\
&&i=-n,\ldots,-1,1,\ldots,n\qquad\qquad\mbox{if }H=U_qso(2n).
\eea
Then the commutation relations (\ref{xixirel}) explicitly amount to
\be 
\ba{ll} 
q\xi^i\xi^j+\xi^j\xi^i=0\qquad\quad &i<j\neq -i,          \\ 
\xi^i\xi^i=0\qquad\quad &i\neq 0,               \\ 
\xi^l\xi^{-l}+\xi^{-l}\xi^{l}=(q-q^{-1}) 
\sum\limits_{i>l }q^{l+1-i}\xi^{-i}\xi^i\qquad\quad &l\ge 1,     \\ 
\xi^0\xi^0=(q^{\frac 12}-q^{- \frac 12}) 
\sum\limits_{i>0} q^{1-i}\xi^{-i}\xi^i. & 
\ea                                        \label{explxixirel} 
\ee
Of course (\ref{explxixirel})$_4$ applies only to 
the cases $H=U_qso(2n\!+\!1)$, and (\ref{explxixirel})$_3$ applies only to 
the cases $H=U_qso(N)$.
The latter relations are equivalent to
the equations  (21) given in Ref. \cite{Fio94}, whence
they can be obtained by an easy rearrangement of terms. 

As already said, as a consequence of (\ref{explxixirel})
$\mbox{dim}(\bigwedge^N)= 1$. Setting e.g. 
\be
\ba{ll}
\gamma_N\,d^N\!x:=\xi^1\xi^{2}...\xi^N\qquad\qquad 
&\mbox{if }H=U_qgl(N)\\[6pt]
\gamma_N\,d^N\!x:=\xi^{-n}\xi^{1-n}...\xi^n\qquad\qquad 
&\mbox{if }H=U_qso(N)
\ea                                                \label{defdNx} 
\ee 
one can introduce
the matrix elements of the $H$-covariant  $\varepsilon$- 
(or completely antisymmetric) tensor  up
to the normalization constant $\gamma_N$ by the relation 
\be
\xi^{i_1}\xi^{i_2}...\xi^{i_N}= d^N\!x\:\varepsilon^{i_1i_2...i_N}. 
\label{defqepsilon}
\ee

The $\varepsilon$-tensors enter the definitions of the
``$q$-determinants'' \cite{FadResTak89,Fio94}, special central 
elements in the Hopf algebras $H'$ dual to $H$, namely the
algebras of functions on the quantum groups. In the appendix
we prove the following proposition, which states a similar 
property for the $q$-determinat of the matrices ${\cal L}^{\pm}$
having as matrix elements the socalled FRT (Faddeev-Reshetikin-Takhtadjan)
generators  \cite{FadResTak89} of $U_qsl(N),U_qso(N)$:
\be 
{\cal L}^{+,}{}_l^a:=\R^{(1)}\rho_l^a(\R^{(2)}),\qquad\qquad 
{\cal L}^{-,}{}_l^a:=\rho_l^a(\R^{-1}{}^{(1)})\R^{-1}{}^{(2)}; \label{frt} 
\ee 
here $\R$ denotes the quasitriangular structure.

\begin{prop} 
\be
\ba{l}
{\cal L}^{+,}{}^{j_N}_{i_N}...{\cal L}^{+,}{}^{j_1}_{i_1} 
\varepsilon^{i_1...i_N}=\varepsilon^{j_1...j_N},\\[6pt] 
{\cal L}^{-,}{}^{j_N}_{i_N}...{\cal L}^{-,}{}^{j_1}_{i_1} 
\varepsilon^{i_1...i_N}=\varepsilon^{j_1...j_N}. 
\ea                                         \label{L-det} 
\ee 
In particular 
\be
\ba{ll}
\det_q{\cal L}^{\pm}:={\cal L}^{\pm,}{}^N_{i_N}...{\cal
L}^{\pm,}{}^{1}_{i_1}  \varepsilon^{i_1...i_N}=\gamma_N
\qquad\qquad &\mbox{if }H=U_qsl(N),\\[6pt] 
\det_q{\cal L}^{\pm}:={\cal L}^{\pm,}{}^n_{i_N}...{\cal L}^{\pm,}{}^{-n}_{i_1} 
\varepsilon^{i_1...i_N}=\gamma_N\qquad\qquad &\mbox{if } H=U_qso(N).
\ea
\ee 
\label{L-det'} 
\end{prop} 

The $U_qso(N)$-covariant matric matrix introduced in (\ref{Pt}) 
coincides with its inverse and
is given by 
\be 
g_{ij}=g^{ij}=q^{\rho_j} \delta_{-i,j}.          \label{defgij} 
\ee 
where 
$(2\rho_j):=(N\!-\!2,N\!-\!4,\ldots,1,0,
-1, \ldots,2\!-\!N)$ 
for $N$ odd,  
$(2\rho_j):=(N\!-\!2,N\!-\!4,\ldots,0,0,\ldots,2\!-\!N)$  
for $N$ even. Introducing the matrix $U$ by 
\be
\ba{ll}
U^i_j=\delta^i_jq^{2i\!-\!N\!-\!1}\qquad\qquad &\mbox{if }H=U_qgl(N)\\
U^i_j:=g^{ik}g_{jk}=\delta^i_jq^{-2\rho_i}\qquad &\mbox{if }H=U_qso(N)
\ea
\ee
(note that $\det U=1$), we can also recall the $q$-cylic property
\cite{Ste96}\footnote{The proof
 given in Ref. \cite{Ste96} applies also to
the $H=U_qgl(N)$ case.}  
 \be
\epsilon^{i_1 ... i_N} = (-1)^{N-1} U^{i_i}_{j_1}  
\epsilon^{i_2 ... i_N j_1}.                       \label{cycl_eps} 
\ee

Let $I:=(i_1,...,i_N)$, and if $I$ is  a permutation of
$J:=(1,...,N)$ denote by $l(I)$ its `length', namely its number of
inversions. The $U_qgl(N)$-covariant deformation of the
$\varepsilon$-tensor \cite{Lyu87,Gur88,FadResTak89} admits the following
compact expression, which closely resembles the undeformed counterpart:
\be
\varepsilon^{i_1i_2...i_N}=\cases{(-q)^{l(I)}\quad \mbox{if $I$ is 
a permutation of $J$},\cr
0\qquad\quad \mbox{otherwise},}
\ee

For the $U_qso(N)$-covariant one  \cite{Fio94} so far no such compact
espression has been found.
In \cite{Fio94,Maj95} several
properties  for general $N$ and
the explicit expression  
for $\varepsilon$ in the cases $N=3,4$
have been determined; we  rewrite them here: for $N=3$,
with normalization $\gamma_3=q^{-1}$
$$
\begin{array}{|c|c|c|c|}
\hline
\varepsilon^{-101}=q^{-1} & \varepsilon^{-110}=-1 &
\varepsilon^{0-11}=-1 & \varepsilon^{01-1}=1 \\ \hline
\varepsilon^{10-1}=-q & \varepsilon^{1-10}=1 & 
\varepsilon^{000}=1/\sqrt{q}-\sqrt{q} & 
\varepsilon^{ijk}=0~~~\mbox{otherwise}, \\
\hline     
\end{array}
$$
and for $N=4$, with normalization $\gamma_4=q^{-2}$
$$
\begin{array}{|c|c|c|c|}
\hline
\varepsilon^{-2-112}=q^{-2} & \varepsilon^{-21-12}=-q^{-2}
&\varepsilon^{-2-121}=-q^{-1} &\varepsilon^{-212-1}=q^{-1} \\  \hline
\varepsilon^{-22-11}=1 & \varepsilon^{-221-1}=-1 &
\varepsilon^{-1-212}=-q^{-1} &  \varepsilon^{-11-22}=1 \\ 
\hline
\varepsilon^{-1-221}=1 & \varepsilon^{-12-21}=-1 &
\varepsilon^{-121-2}=q &  \varepsilon^{-112-2}=-1 \\ 
\hline
\varepsilon^{1-1-22}=-1 & \varepsilon^{1-2-12}=q^{-1} &
\varepsilon^{1-12-2}=q &  \varepsilon^{12-1-2}=-q \\ 
\hline
\varepsilon^{12-2-1}=1 & \varepsilon^{1-22-1}=-1 &
\varepsilon^{2-2-11}=-1 &  \varepsilon^{2-1-21}=q \\ 
\hline
\varepsilon^{21-2-1}=-q & \varepsilon^{2-21-1}=1 &
\varepsilon^{2-11-2}=-q^2 &  \varepsilon^{21-1-2}=q^2 \\ 
\hline
\varepsilon^{-11-11}=k& \varepsilon^{1-11-1}=-k&
\multicolumn{2}{|c|}{\varepsilon^{ijkl}=0~~~\mbox{otherwise}.} \\
\hline
\end{array}
$$
For general $N$ we can at least state the following properties,
which can be easily proved as a consequence of
(\ref{explxixirel}):
\bea
&&\mbox{{\bf Property.} Let } I=\{i_1,...,i_N\}, \: J=\{1,...,N\}.  
                                                 \label{qeps} \\ 
&&\ba{l} 
\mbox{ Then } \varepsilon^{i_1...i_N}=0 
\mbox{ unless all the following conditions are fulfilled: }  \\ 
1.\quad\mbox{ if $N$ is odd, the subset } J_0=\{j\: \:|\: \:i_j=0\} 
\mbox{ has an odd number of}\\ \qquad\mbox{  elements;}\\ 
2.\quad J-J_0 \mbox{ is an union of pairs } \{h,k\} \mbox{ such that }  
i_h=-i_k;\\ 
3.\quad \mbox{the number $\natural_l$ of pairs } \{h,k\} \mbox{ such that
}   i_h=-i_k=l \mbox{ fulfills }\\
\qquad\natural_l\le n\!-\!l\!+\!1;\\ 
4.\quad \mbox{for no }j\in J-(J_0\cup\{N\})\quad i_j=i_{j+1}.
\ea\nonumber
\eea
\be 
\mbox{{\bf Property.} \cite{Fio94}}\qquad\quad
g_{i_1j_1}...g_{i_Nj_N}\varepsilon^{j_N...j_1}=:
\varepsilon_{i_1i_2...i_N}=\varepsilon^{i_N...i_2i_1}.
\ee

\medskip
We now give the relation connecting the antisymmetric projectors and
the $\varepsilon$-tensors. In the appendix we prove
\begin{prop} Let
$(d_0)^{-1}:=\sum\limits_{\{a_i\}}(\varepsilon^{a_1...a_N})^2$. Then
\be
\ba{lll}
{\cal P}^{-,l}{}^{i_1...i_l}_{j_1...j_l}&=&d_l 
U^{k_1}_{j_1}...U^{k_l}_{j_l} \varepsilon^{i_{l\!+\!1}\ldots i_Ni_1\ldots i_l}
\varepsilon^{i_{l\!+\!1}\ldots i_Nk_1\ldots k_l}\\[6pt]
&\stackrel{(\ref{cycl_eps})}{=}& (-)^{l(N\!-\!1)}d_l \varepsilon^{j_1\ldots
j_li_{l\!+\!1}\ldots i_N} 
\varepsilon^{i_{l\!+\!1}\ldots i_Ni_1\ldots i_l} 
\ea \label{Pa_eps_rel0}
\ee
where $d_l$ is defined by
\be
\ba{ll}
d_l := d_0\frac{[N]_q!}{[l]_q![N\!-\!l]_q!} \qquad \qquad
&\mbox{if }H=U_qgl(N),\\[8pt]
d_l := d_0\frac{[N]_q!}{[l]_q![N\!-\!l]_q!} 
\frac{q^{l\!-\!\frac N2}+q^{\frac N2\!-\!l}} 
{q^{-\frac N2}+q^{\frac N2}}
\qquad \qquad &\mbox{if }H=U_qso(N).  
\ea                                             \label{defd_p}  
\ee
Clearly $d_l=d_{N\!-\!l}$, in particular $d_0=d_N$. In the $H=U_qso(N)$ case
this can be also rewritten in the form
\be
{\cal P}^{-,l}{}^{i_1...i_l}_{j_1...j_l}=d_l
\varepsilon_{j_l...j_1}{}^{i_{l\!+\!1}...i_N}
\varepsilon_{i_N...i_{l\!+\!1}}{}^{i_1...i_l}.        \label{Pa_eps_rel}
\ee
\label{Pa-eps}
\end{prop}
By an explicit calculation one finds that for the $\gamma_3,\gamma_4$
given above
\be
d_0^{-1}=\left\{\ba{l} 
[2]_{q^{1/2}}[3]_{q^{1/2}}\quad \:\mbox{ if }N=3\\[8pt]
2([2]_{q^{1/2}})^2[3]_q\quad\mbox{if }N=4
\ea\right.
\ee

\sect{Vielbein bases}
\label{vielb}

The set of $N$ exact forms $\{\xi^i\}$ is a natural basis 
for the $F$-bimodule $\Omega^1$, as well as for the  
the $F\cocross H$-bimodule   $\Omega^1\cocross H$. 
By (\ref{xxirel}), the $\xi^i$ {\it do not} commute with $F$.
We are going to introduce alternative, socalled ``frames'' 
(or ``vielbein''  bases) \cite{DimMad96} which {\it do}, 
revisiting the notion and construction given in 
Ref. \cite{CerFioMad00}.

As shown in \cite{Fio95cmp,ChuZum95}, there exist a 
algebra homomorphism 
\be 
\varphi: \A\cocross H\to \A,              \label{Hom1} 
\ee 
acting as the identity on $\A$ itself, 
\be 
\varphi(a)=a \qquad\qquad a\in \A,     \label{Hom1'} 
\ee 
where $H$ is either Hopf algebra $H=U_qsl(N), U_qso(N)$ and $\A={\cal H}$ 
is the corresponding deformed Heisenberg algebra on $\b{C}_q^N,\b{R}_q^N$.
One can immediately extend $\varphi$ to the central extensions
$H=U_qgl(N),\widetilde{U_qso(N)}$ by setting
\be  
\varphi(q^{-2\eta})=q^N\Lambda^{-2}      \label{varphieta} 
\ee
(adopting the same normalization factor $q^N$ as in \cite{Fio03}), where
the element $\Lambda^{-2}\in {\cal H}$ is defined by \cite{Ogi92}
\bea
&&\Lambda^{-2}\!:=\!1 \!+\!qkx^i\partial_i\equiv 1\!+\!O(h) \qquad\qquad 
  \qquad\qquad\mbox{if } H=U_qsl(N),\\
&&\Lambda^{-2}\!:=\!1 \!+\!qkx^i\partial_i\!+\! 
\frac{q^N\,k^2}{(1\!+\!q^{N-2})^2}r^2\Box\equiv 1\!+\!O(h)\qquad
\mbox{if } H=U_qso(N),\qquad\quad
\eea
(in   \cite{Ogi92,OgiZum92} it was denoted by $\Lambda$). We are
also extending ${\cal H}$ so as to contain
its square root $\Lambda^{-1}$ and inverse square root 
$\Lambda$  as additional 
generators or as formal power series in the deformation parameter 
$h=\ln q$. The latter fulfill the relations 
\be 
\Lambda x^i=q^{-1}x^i\Lambda,\qquad\quad 
\Lambda\partial^i=q\partial^i\Lambda, \qquad\quad 
\Lambda \xi^i=\xi^i\Lambda,  \qquad\label{Lambdaprop} 
\ee 
and the corresponding ones for $\Lambda^{-1}$. 
For real $q$, $\varphi$ is even a ${\star}$-algebra homomorphism. 
Applying $\varphi$ in particular to both sides of (\ref{crossprod}) one 
finds 
\be 
a \,\varphi(g)=\varphi(g_{(1)})\, (a\tl g_{(2)}). \label{crossprod'} 
\ee 

In Ref. \cite{CerFioMad00} to introduce a frame on $\b{R}_q^N$
first auxiliary objects 
\be
\vartheta^i:=q^{-\eta}{\cal L}^{-,}{}^i_l\xi^l\,\in \Omega^1\cocross H
\ee
[with $H=\widetilde{U^-_qso(N)}$, the negative Borel subalgebra of 
$\widetilde{U_qso(N)}$] were introduced, characterized by the property 
to commute with  $F$
\be 
[\vartheta^i,F]=0.                              \label{framecond} 
\ee 
The reader can check (\ref{framecond}) 
by a direct computation that $[\vartheta^i,x^j]=0$.
In \cite{CerFioMad00} we also showed that there exists a
suitable map  $\varphi^-$ 
of the type (\ref{Hom1}-\ref{Hom1'}), with $\A$ a slight extension
of $F$ and  $H=\widetilde{U^-_qso(N)}$.
Replacing $q^{-\eta}{\cal L}^{-,}{}^i_l$ by its $\varphi^-$-image has no
effect on the commutation relation with $x^j$,  see (\ref{crossprod'}),
whence we found that the elements 
$\tilde\theta^i:=\varphi^-(q^{-\eta}{\cal L}^{-,}{}^i_l)\xi^l$  
(in \cite{CerFioMad00} denoted simply as $\{\theta^i\}$)
also fulfilled (\ref{framecond}), and therefore were called 
elements of a ``frame'' (or ``vielbein''), according to the notion
introduced in Ref. 
\cite{DimMad96}.  Now it is also easy to check that 
the {\it same} $\vartheta^i$ also commute with
the derivatives, $[\vartheta^i,\partial^j]=0$.  By the same reasoning,
replacing in the theorems and proofs of  Ref. \cite{CerFioMad00} the map
$\varphi^-$ with   the one $\varphi$,   we arrive at

\begin{prop} 
The sets $\{\vartheta^i\}$ and $\{\theta^i\}$ of 1-forms given by 
\bea 
&&\vartheta^i:=q^{-\eta-\frac N2}{\cal L}^{-,}{}^i_l\xi^l=  
\xi^m q^{1-\eta} \rho_m^j(u_4) {\cal L}^{-,}{}^i_j, 
                                         \label{defvarthetaa}\\ 
&&\theta^i:=\Lambda^{-1}\varphi({\cal L}^{-,}{}^i_l)\xi^l 
=\xi^h\Lambda^{-1}\varphi(S^2{\cal L}^{-,}{}^i_h)   \label{defthetaa'} 
\eea 
($u_4:=\R^{-1}{}^{(1)}S^{-1}\R^{-1}{}^{(2)}$) are resp. ``frame'' bases of 
the ${\cal H}\cocross \widetilde{U_qso(N)}$-bimodule 
${\cal DC}^1\cocross \widetilde{U_qso(N)}$ and of the ${\cal H}$-bimodule 
${\cal DC}^1$, in the sense that 
\be 
[\vartheta^i,{\cal H}]=0\qquad\qquad[\theta^i,{\cal H}]=0.  
                                              \label{framecond'} 
\ee 
They satisfy the same commutation relations as the $\xi^i$, 
\be 
(\1_{N^2}-{\cal P}^t)^{ij}_{hk}\vartheta^h\vartheta^k=0 \qquad\qquad  
(\1_{N^2}-{\cal P}^t)^{ij}_{hk}\theta^h\theta^k=0. \label{ththrel} 
\ee 
Finally, they form a $N$-plet under the action 
of $U^{op}_qso(N)$ (i.e. $U_qso(N)$ endowed with the opposite 
coproduct). 
\end{prop} 
 \noindent  
We just give the proof of the second equality in (\ref{defthetaa'}), which
was not given in   \cite{CerFioMad00}. Recalling the coproduct
$\Delta({\cal L}^{-,}{}^i_k)={\cal L}^{-,}{}^i_h\otimes {\cal L}^{-,}{}^h_k$ 
of the FRT generators we find
\bea  
\theta^i &=& 
\varphi(S{\cal L}^{-,}{}^h_k\,S^2{\cal L}^{-,}{}^i_h)\theta^k 
\stackrel{(\ref{framecond'})_2}{=} 
\varphi(S{\cal L}^{-,}{}^h_k)\theta^k\varphi(S^2{\cal L}^{-,}{}^i_h)\nn 
&\stackrel{(\ref{defthetaa'})}{=}& 
\xi^h\Lambda^{-1}\varphi(S^2{\cal L}^{-,}{}^i_h) 
=\xi^h\Lambda^{-1}U_h^j\varphi({\cal L}^{-,}{}^k_j)U ^{-1}{}^i_k.\nonumber 
\eea  

\medskip
An analogous proposition for objects $\hat\vartheta^a,\hat\theta^a$ 
obtained by replacing $\xi^i$ by $\hat\xi^i$, ${\cal L}^{-,}{}^i_h$ by 
${\cal L}^{+,}{}^i_h$ and $\eta$ by $-\eta$ holds.
 In Ref. \cite{Fio03}
we have shown that the frame basis elements $\hat\theta^a$ transform
exactly as the coordinates $x^a$ under the $\star$-structure
chraterizing real $q$ (namely can be made real by a suitable $\b{C}$-linear
transformation).

\medskip
Explicit expressions for the images $\varphi({\cal L}^{-,}{}^h_k)$
for the $H$'s with the lowest $N$'s, $H=U_qsl(2),U_qso(3)$, can be
found e.g. in Ref. \cite{FioSteWes03}, section 4.

As the commutation relations (\ref{ththrel}) among the $\theta^h$ are exactly of
the same form of the ones (\ref{xixirel}) among the $\xi^i$, we immediately
find  that also the space $\bigwedge^N_{\theta}$ of monomials in  
$\theta^i$ of degree $N$ has dimension 1. Moreover,
\be
\theta^{i_1}\theta^{i_2}...\theta^{i_N}=  
dV\:\varepsilon^{i_1i_2...i_N}. 
\ee
where $dV\in\mbox{$\bigwedge^N_{\theta}$}$ is defined replacing in 
(\ref{defdNx}) $d^N\!x$ with $dV$ and $\xi^i$ with $\theta^i$.

\begin{prop} 
\label{volform} 
The ``volume form'' $dV$ is central in ${\cal DC}^*$ and equal to 
\be 
dV=d^N\!x\Lambda^{-N}. 
\ee 
\end{prop} 
{\bf Proof:} With the definition of $dV$ adopted for $H=\widetilde{U_qso(N)}$ 
(the case $H=U_qgl(N)$ is completely analogous)
\bea 
dV&\stackrel{(\ref{defthetaa'})}{=}& 
\theta^{-n}...\theta^{n\!-\!1} 
\Lambda^{-1}\varphi({\cal L}^{-,}{}^n_{i_N})\xi^{i_N} 
\stackrel{(\ref{framecond'})_2}{=} 
\Lambda^{-1}\varphi({\cal L}^{-,}{}^n_{i_N}) 
\theta^{-n}...\theta^{n\!-\!1}\xi^{i_N}\nn 
&=&...=\Lambda^{-1}\varphi({\cal L}^{-,}{}^n_{i_N})... 
\Lambda^{-1}\varphi({\cal L}^{-,}{}^{-n}_{i_1})\xi^{i_1}...\xi^{i_N}\nn 
&\stackrel{(\ref{defqepsilon})}{=}&\Lambda^{-N}\varphi({\cal L}^{-,}{}^n_{i_N} 
...{\cal L}^{-,}{}^{-n}_{i_1})\varepsilon^{i_1...i_N}d^N\!x 
\stackrel{(\ref{L-det}),(\ref{varphieta})}{=}\Lambda^{-N}d^N\!x 
\qquad\qquad\Box\nonumber 
\eea

The reader might wonder about the usefulness of the generalized notion of
vielbein introduced in this section: generally
speaking the  differential forms $\omega_p\in\Omega^p$ and the functions
$f\in F$ have a geometrical or physical significance, so since $\theta^a$ are
in ${\cal DC}^1$ rather than in $\Omega^1$, the components of $\omega_p$ in the
vielbein basis are in ${\cal H}$ rather  than in $F$. The point is that, as we
have shown in Ref. \cite{Fio03}, the difference between these components is
irrelevant when evaluating functionals on $\Omega^N$, scalar products in
$\Omega^p$, etc.  by means of integration, provided Stokes' theorem applies.

\sect{Hodge map and Laplacian on $\b{R}_q^N$} 
\label{hodge} 

Having at one's disposal also the $U_qso(N)$-covariant metric matrix $g_{ij}$,
a ($U_qso(N)$-covariant) Hodge map   $*:\bigwedge^p\to\bigwedge^{N-p}$ 
acting on  {\it exterior} forms on $\b{R}_q^N$ was introduced (leaving some
ambiguities) in \cite{fiothesis,Maj95} using both $g_{ij}$ and the $q$-epsilon
tensor, in analogy with the undeformed theory. As we are going to see,  one
has to fix the ambiguities to make $*$ involutive and moreover add in the
definition a suitable power of $\Lambda$ in order to define a Hodge map  on
{\it differential} forms. It is  more convenient to start giving the
definition of the Hodge map  in the frame basis:   

\begin{prop} 
For $H=U_qso(N)$
and any $p=0,1,...,N$ one can define a ${\cal H}$-bilinear map 
\be 
*:{\cal DC}^p\to{\cal DC}^{N-p} 
\ee 
the ``Hodge map'', such that ${}^*\1=dV$ and on each ${\cal DC}^p$, and
therefore on the whole ${\cal DC}^*$,
\be 
*^2\equiv *\circ *=\id       \label{involution} 
\ee 
by setting on the monomials in the $\theta^a$ 
\be 
{}^*(\theta^{a_1}\theta^{a_2}...\theta^{a_p}) 
=\,c_p\,\theta^{a_{p+1}}...\theta^{a_N} 
\varepsilon_{a_N...a_{p+1}}{}^{a_1...a_p};           \label{defHodge1} 
\ee 
the normalization constants  $c_p$ are 
constrained by the conditions 
\be 
c_p\,c_{N\!-\!p}=d_p.                  \label{condc_p}
\ee  
\end{prop}

The most convenient choice for the $c_p$ will be given below.
${\cal H}$-bilinearity implies in particular 
 
\be 
{}^*(a\,\omega_p\, b)=a\,{}^*\omega_p\, b\qquad\qquad 
\forall \,a,b \in {\cal H}, \quad \omega_p\in{\cal DC}^p;  \label{Hbil} 
\ee 
i.e. applying Hodge and multiplying by ``functions 
or differential operators'' 
are commuting operations, in other words 
a differential form $\omega_p$ and its Hodge map
image have the same commutation 
relations with $x^i,\partial^j$. That this is true is evident in the  
frame basis, because of (\ref{framecond'}).  
Relation (\ref{involution}) easily follows from  (\ref{Pa_eps_rel}).
The fixed positive sign at the rhs of (\ref{involution})
[cumbersome when compared with the more familiar $(-1)^{p(N-p)}$] 
is the sign of $d_0$ and is due to the non-standard 
ordering of the indices in (\ref{defHodge1}). The latter in 
turn is the only correct one: had we used a different 
order, at the rhs of (\ref{involution}) tensor products 
of the matrices $U^{\pm 1}$, instead of the identity, would have 
appeared, because of property (\ref{cycl_eps}). 
 
Using the ${\cal H}$-bilinearity of $*$ in the appendix we prove 
 
\begin{prop} 
In terms of the basis of differentials (\ref{defHodge1}) takes the form 
\be 
{}^*(\xi^{i_1}...\xi^{i_p})=q^{-N(p-N/2)} 
c_p\xi^{i_{p+1}}...\xi^{i_N} 
\varepsilon_{i_N...i_{p+1}}{}^{i_1...i_p}\Lambda^{2p-N}. 
 \label{defHodge2} 
\ee 
\label{Hodge2} 
\end{prop} 
 
This differs from the (incomplete) definition of Hodge map on exterior forms  
given in \cite{fiothesis,Maj95} by the presence 
of $\Lambda$-powers (needed for the ${\cal H}$-bilinearity), 
by the already noted crucial different indices order and by the explicit 
determination of the 
coefficients $c_p$. From the above formulae and the 
commutation relations 
 (\ref{xxirel}),  (\ref{Lambdaprop}) it is evident that 
by restricting the domain of $*$ to the unital 
subalgebra $\widetilde{\Omega}^*\subset {\cal DC}^*$ generated by 
$x^i,\xi^j,\Lambda^{\pm 1}$ one obtains a  
$\widetilde{F}$-bilinear map 
\be 
*:\widetilde{\Omega}^p\to\widetilde{\Omega}^{N-p} 
\ee 
fulfilling again (\ref{involution})  
[just take $a,b\in\widetilde{F}$ in (\ref{Hbil})]; here  
$\widetilde{F}$ denotes the unital subalgebra 
generated by $x^i,\Lambda^{\pm 1}$.  
This restriction is what is the notion closest to the conventional notion
of a Hodge map on $\b{R}_q^N$: as a matter of fact, there is no
$F$-bilinear restriction of $*$ to $\Omega^*$.
 
{}From the bilinearity of the Hodge map and the explicit 
$U_qso(N)$-covariant form of (\ref{defHodge2}) it immediately follows 
 
\begin{prop} 
The Hodge map is $U_qso(N)$-covariant, i.e. commutes with the 
$U_qso(N)$-action: 
\be 
({}^*\omega_p)\tl g= {}^*(\omega_p\tl g)\qquad\qquad 
\forall g\in U_qso(N) 
\ee 
This is true also for its restriction to the 
subalgebra $\widetilde{\Omega}^*\subset {\cal DC}^*$. 
\end{prop} 
{\bf Remark:} But $*$ is not $\widetilde{U_qso(N)}$-covariant. 
This is due to the fact that $\eta$ has a nontrivial action on 
each $\xi^i$, and $*$ changes the degree of a monomial in the 
$\xi^i$'s. 

As in commutative geometry we introduce the exterior 
coderivative by
\be
\delta:=-{}^*d{}^*.
\ee
In Ref. \cite{Fio03} we show that (at least for 
positive $q$) $\delta$ can be seen as the hermitean conjugate
of $d$ acting on $\Omega^*$ endowed with a suitable
scalar product.
The residual freedom  left by (\ref{condc_p}) in choosing
the $c_p$ is eliminated by requiring that the
differential operator $d\,\delta+\delta\,d$ is a scalar proportional to 
$\partial\cdot\partial$, as in the commutative geometry case.
In the appendix we prove the following proposition: 
 
\begin{prop} 
The ``Laplacian'' $\Delta:=d\,\delta+\delta\,d$ reduces on
all ${\cal DC}^*$, and in particular on $\Omega^*$, to
\be 
\Delta =-q^2\partial\cdot\partial\Lambda^2 =-q^{-N} 
\hat\partial\cdot\hat\partial,            \label{Laplacian}     
\ee
provided we choose
\be
c_p=\frac 1{[N\!-\!p]_q!}
\prod\limits_{l=p}^{N\!-\!1}
\frac{q^{l\!-\!\frac N2}\!+\!q^{\frac N2\!-\!l}}
{q^{1\!-\!\frac N2}\!+\!q^{\frac N2\!-\!1}}.   \label{c_p=}
\ee 
\label{lapl}
\end{prop}

\app{Appendix}

We begin this appendix by recalling few basic properties about the universal
$R$-matrix,  or quasitriangular structure \cite{Dri85-86}, $\R$ 
of the quantum groups $\uqg$, while fixing our conventions. 
$\R$ intertwines between  $\Delta$ and opposite coproduct 
$\Delta^{op}$, and so does also $\R^{-1}_{21}$: 
\be 
\ba{l} 
\R( g_{(1)}\otimes g_{(2)})=(g_{(2)}\otimes g_{(1)})\R,  \\ 
\R^{-1}_{21}(g_{(1)}\otimes g_{(2)})=(g_{(2)}\otimes g_{(1)})\R^{-1}_{21}.   
\ea 
\label{inter} 
\ee 
It fulfills 
\bea 
&&(\Delta \otimes \mbox{id})\R=\R_{13}\R_{23}, \label{delta1}\\ 
&&(\mbox{id}\otimes\Delta  )\R=\R_{13}\R_{12}, \label{delta2}\\ 
&&(S \otimes \mbox{id})\R=\R^{-1}= (\mbox{id}\otimes S^{-1})\R, 
                                               \label{SR} \\ 
&& S^{-1}(g) = u^{-1}S(g) u                  \label{interantip} 
\eea 
where $u$, which is defined up to an invertible central factor, can be taken 
e.g. as the $u=u_1$ with
\be
u_1:= (S\R^{(2)}) \R^{(1)}, \label{defu1}
\ee
{From} (\ref{inter}-\ref{delta2}) it follows the 
universal Yang-Baxter relation 
\be 
\R_{12}\R_{13}\R_{23}=\R_{23}\R_{13}\R_{12}.          \label{YBEQ1} 
\ee 
The braid matrix $\hat R$ \cite{FadResTak89} is related to the quasitriangular
structure $\R$ by $\hat R^{ij}_{hk}\equiv
R^{ji}_{hk}:=q^{\sigma_N}(\rho^j_h\otimes\rho^i_k)\R
= q^{\sigma_N}\rho^j_h\Big( {\cal L}^{+,}{}^i_k\Big)$, where $\sigma_N=1/N$
for $H=U_qgl(N)$ and $\sigma_N=0$ for $H=U_qso(N)$. With the indices'
convention described in section \ref{eps} $\hat R$ is
given by 
\be
\hat R = q \sum_i e^i_i \otimes e^i_i +
\sum_{\scriptstyle i \neq j} e^j_i \otimes e^i_j
+k \sum_{i<j} e^i_i \otimes e^j_j             \label{defRslN}
\ee
when $H=U_qsl(N)$, and by
\bea
\hat R&=&q \sum_{i \neq 0} e^i_i \otimes e^i_i +
\sum_{\stackrel{\scriptstyle i \neq j,-j} 
{\mbox{ or } i=j=0}} e^j_i \otimes e^i_j+ q^{-1} 
\sum_{i \neq 0} e^{-i}_i
\otimes e^i_{-i} 
\label{defRsoN} \\
&&+k (\sum_{i<j} e^i_i \otimes e^j_j- 
\sum_{i<j} q^{-\rho_i+\rho_j} 
e^{-j}_i \otimes e^j_{-i}) \nonumber
\eea
when $H=U_qso(N)$.
Here $e^i_j$ is the $N \times N$ matrix with all elements
equal to zero except for a $1$ in the $i$th column and $j$th row.,
and $k:=q-q^{-1}$.

In the $H=U_qso(N)$ case, using (\ref{gRrel}), (\ref{Pt}) it is not difficult
to show the following formulae 
\be
\ba{l}
{\cal P}^t_{12}\hat R^{\pm 1}_{23}=
Q_N{\cal P}^t_{12}{\cal P}^t_{23}\hat R^{\mp 1}_{12},
\qquad\quad
\hat R^{\pm 1}_{23}{\cal P}^t_{12}=
Q_N\hat R^{\mp 1}_{12}{\cal P}^t_{23}{\cal P}^t_{12},\\
{\cal P}^t_{23}\hat R^{\pm 1}_{12}=
Q_N{\cal P}^t_{23}{\cal P}^t_{12}\hat R^{\mp 1}_{23},
\qquad\quad
\hat R^{\pm 1}_{12}{\cal P}^t_{23}=
Q_N\hat R^{\mp 1}_{23}{\cal P}^t_{12}{\cal P}^t_{23},
\ea
\ee
which are written in matrix-tensor notation in order to let us
do many proofs avoiding indices $i,j$ etc.
 Moreover, 
\be
\rho^a_b(Sh)=g^{ad}\rho^c_d(h)g_{cb},
\qquad\qquad\Rightarrow\qquad\qquad S\c{L}^{\mp,}{}^j_i=
g_{ih}\c{L}^{\mp,}{}^h_k g^{kj}.            \label{antipode}
\ee

\subsubsection*{Proof of Proposition \ref{symmetrizers} }

One can determine the projectors ${\cal P}^{\pm,l}$
iteratively. We adopt the Ansatz (\ref{ansatz1})
with $M^{\pm}=f^{\pm}(\hat R)$ a matrix to be determined.
The most general one is
\be
\ba{ll}
M^{\pm,l\!+\!1}=\alpha^{\pm}_{l_{N^2}\!+\!1}\left(\1+\beta^{\pm}_{l\!+\!1}\hat
R\right) \qquad\qquad\qquad &\mbox{ if } H=U_qsl(N)\qquad \\[8pt]  
M^{\pm,l\!+\!1}=\alpha^{\pm}_{l_{N^2}\!+\!1}\left(\1+\beta^{\pm}_{l\!+\!1}\hat R
+\gamma^{\pm}_{l\!+\!1}{\cal P}^t\right)\qquad\quad &\mbox{ if }
H=U_qso(N)\qquad \ea                                         \label{m2}
\ee
We first determine the coefficients
$\beta^{\pm}_{l\!+\!1},\gamma^{\pm}_{l\!+\!1}$ by imposing the conditions
(\ref{Plproj1}). By the recursive assumption, only the condition with $m=l$ is
not fulfilled automatically and must be imposed by hand. Actually, it suffices
to impose just (\ref{Plproj1})$_1$, due
to the symmetry of the Ansatz (\ref{ansatz1}) and of the matrices 
${\cal P}^{\pi}$ under transposition. Setting
\bea
&&{\cal P}^{'\mp}_{l(l\!+\!1)}:=
{\cal P}^{\mp}_{l(l\!+\!1)}\qquad\qquad\qquad\quad\mbox{ if } H=U_qsl(N),\nn
&&{\cal P}^{'\mp}_{l(l\!+\!1)}:=
{\cal P}^{\mp}_{l(l\!+\!1)}+{\cal P}^t_{l(l\!+\!1)}\qquad\qquad \mbox{ if }
H=U_qso(N), \nonumber
\eea
this amounts to
\bea
0&\stackrel{!}{=}&{\cal
P}^{\pm,l\!+\!1}{\cal P}^{'\mp}_{l(l\!+\!1)}
\stackrel{(\ref{ansatz1})}{=} {\cal P}^{\pm,l}_{1\ldots
l}M^{\pm}_{l(l\!+\!1)}{\cal P}^{\pm,l\!-\!1}_{1...(l\!-\!1)}
M^{\pm}_{(l\!-\!1)l} {\cal P}^{\pm,l\!-\!1}_{1...(l\!-\!1)} {\cal
P}^{'\mp}_{l(l\!+\!1)}\nn   &=& {\cal P}^{\pm,l}_{1\ldots
l}M^{\pm}_{l(l\!+\!1)}M^{\pm}_{(l\!-\!1)l}  {\cal P}^{\pm,l\!-\!1}_{1\ldots
(l\!-\!1)} {\cal P}^{'\mp}_{l(l\!+\!1)}. \label{cond} \eea

In the $H=U_qsl(N)$ case (\ref{cond}) becomes
\bea
0&\propto & {\cal P}^{\pm,l}_{1\ldots
l}\left[\1_{N^{l\!+\!1}}\!+\!\beta^{\pm}_{l\!+\!1}\hat R_{l(l\!+\!1)}\!+\!
\beta^{\pm}_l\hat R_{(l\!-\!1)l}\!+\! \beta^{\pm}_l\beta^{\pm}_{l\!+\!1}\hat R_{l(l\!+\!1)}
\hat R_{(l\!-\!1)l} \right] {\cal P}^{\pm,l\!-\!1}_{1\ldots (l\!-\!1)} {\cal
P}^{'\mp}_{l(l\!+\!1)} \nn 
&=& {\cal P}^{\pm,l}_{1\ldots l}\left[1\!\mp\!
q^{\mp 1}\beta^{\pm}_{l\!+\!1} \!\pm\!
q^{\pm 1}\beta^{\pm}_l\right]
{\cal P}^{\pm,l\!-\!1}_{1\ldots (l\!-\!1)} {\cal P}^{'\mp}_{l(l\!+\!1)}
\label{blu} 
\eea
where we have used the braid equation (\ref{braid2}) to see that
the term proportional to $\beta^{\pm}_l\beta^{\pm}_{l\!+\!1}$ vanishes,
and the relations 
$$
\hat R {\cal P}^{'\mp}=\mp q^{\mp 1}\1_{N^2}{\cal P}^{'\mp},
\qquad\qquad {\cal P}^{\pm,l}_{1\ldots l}\hat R_{(l\!-\!1)l}
=\pm q^{\pm 1}{\cal P}^{\pm,l}_{1\ldots l}.
$$
The condition that the square bracket in (\ref{blu}) vanishes
is recursively solved,  starting from $l=1$ with initial input
$\beta^{\pm}_1=0$ (since ${\cal P}^{\pm,1}=\1_N$), by
\be
\beta^{\pm}_{l\!+\!1}=\pm q^{\pm 1}\, l_{q^{\pm 2}}. \label{beta}
\ee
[for $l=2$ this gives back (\ref{projectorpm1})].

In the $H=U_qso(N)$ case (\ref{cond}) becomes
\bea
0&\propto & {\cal P}^{\pm,l}_{1\ldots l}\left[\1_{N^{l\!+\!1}}\!+\!\beta^{\pm}_{l\!+\!1}\hat
R_{l(l\!+\!1)} \!+\!\gamma^{\pm}_{l\!+\!1}{\cal P}^t_{l(l\!+\!1)}\!+\!
\beta^{\pm}_l\hat R_{(l\!-\!1)l}\!+\!\gamma^{\pm}_l{\cal
P}^t_{(l\!-\!1)l}\right.\nn
&& \!+\! \beta^{\pm}_l\beta^{\pm}_{l\!+\!1}\hat R_{l(l\!+\!1)}
\hat R_{(l\!-\!1)l}\!+\!\gamma^{\pm}_l\gamma^{\pm}_{l\!+\!1}
{\cal P}^t_{l(l\!+\!1)}{\cal
P}^t_{(l\!-\!1)l}\!+\!\beta^{\pm}_l\gamma^{\pm}_{l\!+\!1} {\cal
P}^t_{l(l\!+\!1)}\hat R_{(l\!-\!1)l} \nn
&&\left. \!+\!\gamma^{\pm}_l\beta^{\pm}_{l\!+\!1}
\hat R_{l(l\!+\!1)}{\cal P}^t_{(l\!-\!1)l} \right]
{\cal P}^{\pm,l\!-\!1}_{1\ldots (l\!-\!1)} {\cal P}^{'\mp}_{l(l\!+\!1)} \nn
&=& {\cal P}^{\pm,l}_{1\ldots l}\left\{\1_{N^{l\!+\!1}}\!+\!\beta^{\pm}_{l\!+\!1}
\left[\mp q^{\mp 1}\1_{N^{l\!+\!1}}
\!+\!(q^{1\!-\!N}\pm q^{\mp 1}){\cal P}^t_{l(l\!+\!1)}\right]
\!+\!\gamma^{\pm}_{l\!+\!1}{\cal P}^t_{l(l\!+\!1)}\!\pm\!
q^{\pm 1}\beta^{\pm}_l\1_{N^{l\!+\!1}}\right.\nn
&&\!+\!\gamma^{\pm}_l\gamma^{\pm}_{l\!+\!1}
{\cal P}^t_{l(l\!+\!1)}{\cal
P}^t_{(l\!-\!1)l}\!+\!\beta^{\pm}_l\gamma^{\pm}_{l\!+\!1} Q_N {\cal
P}^t_{l(l\!+\!1)}{\cal P}^t_{(l\!-\!1)l}\left[\mp q^{\pm 1}\1_{N^{l\!+\!1}}
\!+\!\left(q^{N\!-\!1}\right.\right.
\nn
&& \left. \left.\left.  \pm q^{\pm 1}\right){\cal P}^t_{l(l\!+\!1)}\right]
\!\pm\! q^{\mp 1}\gamma^{\pm}_l\beta^{\pm}_{l\!+\!1}Q_N
{\cal P}^t_{l(l\!+\!1)}{\cal P}^t_{(l\!-\!1)l}
\right\}
{\cal P}^{\pm,l\!-\!1}_{1\ldots (l\!-\!1)}{\cal P}^{'\mp}_{l(l\!+\!1)} \nn
&=& {\cal P}^{\pm,l}_{1\ldots
l}\left\{\1_{N^{l\!+\!1}}\left[1\!\mp\!q^{\mp 1}\beta^{\pm}_{l\!+\!1}\!\pm\!
q^{\pm 1}\beta^{\pm}_l\right]\!+\!{\cal P}^t_{l(l\!+\!1)} \left[
\beta^{\pm}_{l\!+\!1}(q^{1\!-\!N}\pm q^{\mp 1})\!+\!
\gamma^{\pm}_{l\!+\!1}\right.\right.\nn
&&\left.\!+\!\beta^{\pm}_l\gamma^{\pm}_{l\!+\!1}
\frac{q^{N\!-\!1}\!\pm \!q^{\pm 1}}{Q_N} \right] 
\!+\! {\cal P}^t_{l(l\!+\!1)}{\cal P}^t_{(l\!-\!1)l}
\left[\gamma^{\pm}_l\gamma^{\pm}_{l\!+\!1}
\!\mp\!\beta^{\pm}_l\gamma^{\pm}_{l\!+\!1} Q_N  q^{\pm 1}
\right.\nn
&&\left.\left.
\!\pm\!q^{\mp 1}\gamma^{\pm}_l\beta^{\pm}_{l\!+\!1} Q_N\right]
\right\} {\cal P}^{\pm,l\!-\!1}_{1\ldots (l\!-\!1)} {\cal P}^{'\mp}_{l(l\!+\!1)}.
\nonumber
\eea
where we have used the braid equation (\ref{braid2}) to see that
the term proportional to $\beta^{\pm}_l\beta^{\pm}_{l\!+\!1}$ vanishes,
and the relations
\bea
&& \hat R {\cal P}^{'\mp}=\left[\mp q^{\mp 1}\1_{N^2}
\!+\!(q^{1\!-\!N}\pm q^{\mp 1}){\cal P}^t\right]{\cal P}^{'\mp}\nn[8pt]
&& {\cal P}^t_{l(l\!+\!1)}\hat R_{(l\!-\!1)l}{\cal P}^{'\mp}_{l(l\!+\!1)}=
Q_N {\cal P}^t_{l(l\!+\!1)}{\cal P}^t_{(l\!-\!1)l}\hat
R^{-1}_{l(l\!+\!1)}{\cal P}^{'\mp}_{l(l\!+\!1)} \nn
&& \qquad\qquad \qquad\qquad = Q_N {\cal
P}^t_{l(l\!+\!1)}{\cal P}^t_{(l\!-\!1)l}\!\left[\mp q^{\pm 1}\1_{N^{l\!+\!1}}
\!+\!\left(q^{N\!-\!1}\!\pm\! q^{\pm 1}\right)\!{\cal P}^t_{l(l\!+\!1)}\right]\!{\cal
P}^{'\mp}_{l(l\!+\!1)} \nn[8pt] 
&& {\cal P}^{\pm,l}_{1\ldots l}\hat
R_{l(l\!+\!1)}{\cal P}^t_{(l\!-\!1)l}\!=\!Q_N {\cal P}^{\pm,l}_{1\ldots l}\hat
R^{-1}_{(l\!-\!1)l}{\cal P}^t_{l(l\!+\!1)}{\cal P}^t_{(l\!-\!1)l} \!=\!\pm\! q^{\mp
1}Q_N{\cal P}^{\pm,l}_{1\ldots l}{\cal P}^t_{l(l\!+\!1)}\!{\cal
P}^t_{(l\!-\!1)l}\nn 
&& {\cal P}^t_{l(l\!+\!1)}{\cal P}^t_{(l\!-\!1)l}{\cal P}^t_{l(l\!+\!1)}
=\frac 1{Q_N^2}{\cal P}^t_{l(l\!+\!1)}.
\nonumber 
\eea
The conditions that the three square brackets vanish
\bea
&& 1\!\mp\!q^{\mp 1}\beta^{\pm}_{l\!+\!1}\!\pm\!q^{\pm 1}\beta^{\pm}_l=0,\nn
&& \beta^{\pm}_{l\!+\!1}(q^{1\!-\!N}\pm q^{\mp 1})\!+\!
\gamma^{\pm}_{l\!+\!1}\!+\!\beta^{\pm}_l\gamma^{\pm}_{l\!+\!1}
\frac{q^{N\!-\!1}\!\pm \!q^{\pm 1}}{Q_N}=0, \nn
&& \gamma^{\pm}_l\gamma^{\pm}_{l\!+\!1}
\!\mp\!q^{\pm 1}\beta^{\pm}_l\gamma^{\pm}_{l\!+\!1} Q_N  
\!\pm\!q^{\mp 1}\gamma^{\pm}_l\beta^{\pm}_{l\!+\!1} Q_N=0,\nonumber
\eea
are recursively solved,  starting from $l=1$ with initial input
$\beta^{\pm}_1=0=\gamma^{\pm}_1$ (since ${\cal P}^{\pm,1}=\1_N$), 
again by (\ref{beta}) and by
\be
\gamma^+_{l\!+\!1}=\frac{\left(q^N\!-\!1\right)\!\left(1\!+\!q^{
2\!-\!N}\right)} {1\!-\!q^{N\!+\!2l\!- \!2}}\,l_{q^{ 2}}
\qquad\quad
\gamma^-_{l\!+\!1}=\frac{\left(q^{-N}\!-\!1\right)\!\left(1\!+\!q^{
N\!-\!2}\right)} {1\!-\!q^{N\!-\!2l}}\,l_{q^{ -2}}
\ee
[for $l=2$ this gives back (\ref{projectorpm2})].

We determine the coefficient $\alpha^{\pm}_{l\!+\!1}$ by imposing the
condition (\ref{Plproj2}). For both $H=U_qsl(N),U_qso(N)$ this gives
\bea
0&\stackrel{!}{=}&{\cal P}^{\pm,l\!+\!1}\left({\cal
P}^{\pm,l\!+\!1}-\1_{N^{l\!+\!1}}\right) \stackrel{(\ref{ansatz1})}{=}{\cal
P}^{\pm,l\!+\!1}\left( {\cal P}^{\pm,l}_{1\ldots l}M^{\pm}_{l(l\!+\!1)}
{\cal P}^{\pm,l}_{1\ldots l}-\1_{N^{l\!+\!1}}\right)\nn
&= & {\cal P}^{\pm,(l\!+\!1)}\left[
\alpha^{\pm}_{l\!+\!1}\left(1\pm q^{\pm 1}\beta^{\pm}_{l\!+\!1}\right)-1\right];
\nonumber  
\eea
in the last equality we have used
(\ref{Plproj1}), (\ref{m2}), (\ref{projectorR}). The condition that the square
bracket  vanishes is recursively solved,  starting from $l=0$ with initial
input $\alpha^{\pm}_0=1$, by
$$
\alpha^{\pm}_{l\!+\!1}=\frac 1{(l\!+\!1)_{q^{\pm 2}}}.
$$
By using (\ref{defnum}) we give at the form (\ref{Ml}) for $M^{\pm,l\!+\!1}$.

To check that (\ref{Plproj3}) is satisfied we just note that the dimension of
each projector is an integer, that it is the required one for $q=1$
(since in this limit the projector reduces to its undeformed counterpart),
and therefore it is also for any generic $q$, by continuity in $q$.

\subsubsection*{Proof of Proposition \ref{L-det'}}

Being $H$-invariant, the element   
$\xi^{i_1}\xi^{i_2}...\xi^{i_N}\in \bigwedge^N$ commutes with   
all $H$ (within ${\cal DC}^*\cocross H$). Therefore  
\bea  
\xi^{i_1}\ldots\xi^{i_N} \! 
&\! \stackrel{(\ref{defu1})}{=}\! & u_1^{-1}  
(S\R^{(2)}) \,\xi^{i_1}...\xi^{i_N}\R^{(1)}\nn  
&\stackrel{(\ref{crossprod})}{=}& u_1^{-1}(S\R^{(2)}) \R^{(1)}_{(1)}  
(\xi^{i_1}\tl \R^{(1)}_{(2)})...(\xi^{i_N} \tl \R^{(1)}_{(N+1)})\nn  
&\stackrel{(\ref{fundrep})}{=}& u_1^{-1}(S\R^{(2)}) \R^{(1)}_{(1)}  
\rho^{i_1}_{j_1}(\R^{(1)}_{(2)}) ... \rho^{i_N}_{j_N}(\R^{(1)}_{(N+1)})  
\xi^{j_1}...\xi^{j_N}\nn  
&\stackrel{(\ref{delta1})}{=}\! &\!  u_1^{-1} \!\left[\!S\R^{(2)}_N\!\right]\!  ... 
\!\left[\!S\R^{(2)}_1\!\right]\!\! \left[\!S\R^{(2)}\!\right] \!\R^{(1)}  
\!
\rho^{i_1}_{j_1}\!\left[\!\R^{(1)}_1\!\right]\! ... \!\rho^{i_N}_{j_N}\!\left[\!\R^{(1)}_N\!\right] \!
\xi^{j_1}...\xi^{j_N}\nn  
&\stackrel{(\ref{interantip})}{=}& (S^{-1} \R^{(2)}_N)  ... 
(S^{-1} \R^{(2)}_1) \rho^{i_1}_{j_1}(\R^{(1)}_1) ... \rho^{i_N}_{j_N}(\R^{(1)}_N)  
\xi^{j_1}...\xi^{j_N}\nn  
&\stackrel{(\ref{SR})}{=}& \R^{-1(2)}_N  ... 
\R^{-1(2)}_1 \rho^{i_1}_{j_1}(\R^{-1(1)}_1)   ... 
\rho^{i_N}_{j_N}(\R^{-1(1)}_N)\xi^{j_1}...\xi^{j_N}\nn  
&\stackrel{(\ref{frt})}{=}& {\cal L}^{-,}{}^{i_N}_{j_N} ... 
{\cal L}^{-,}{}^{i_1}_{j_1} \xi^{j_1}...\xi^{j_N},\nonumber  
\eea  
where $\R_1,...,\R_N$ just denote $N$ different copies of $\R$;  
factoring out $d^N\!x$ [see Eq. (\ref{defqepsilon})]  
the claim follows. 

\subsubsection*{Proof of Proposition \ref{Pa-eps}}

We start by recalling the relations  (which can be easily checked
using the explicit definition of $\hat R,U,{\cal P}^t$ given above)
\bea 
&&\mbox{tr}_2\left(U_2\hat R_{12}\right)=q^N\1_N 
\qquad \mbox{tr}(U)=[N]_q \qquad\qquad  \mbox{if }H=U_qgl(N)\nn 
&&\mbox{tr}_2\!\left(\!U_2\!\hat R_{12}\!\right)\!=\!q^{N\!-\!1}\1_N 
\qquad \mbox{tr}_2\!\left(\!U_2\!{\cal P}^t_{12}\!\right)\!=\!\frac{\1_N}{Q_N} 
\qquad \mbox{tr}\!(\!U\!)\!=\!Q_N,  \quad  \mbox{if }H\!=\!U_qso(N)\nonumber 
\eea 
where $U^{\pm 1}_2\equiv\1_N\otimes U^{\pm 1}$ and $\mbox{tr}_2$ denotes matrix
trace on the second factor in the tensor product $\b{C}^N\otimes \b{C}^N$; this
implies
\be 
\mbox{tr}_2\left(U_2M^{-,l}_{12}\right)=b_l \1_N,      \label{tool1} 
\ee 
where 
\bea 
b_l &=& \frac{[N\!-\!l\!+\!1]_q}{[l]_q}\qquad \qquad 
\mbox{if }H=U_qgl(N)\\[8pt] 
b_l &=&\frac 1{[l]_q}\left[q^{l\!-\!1}Q_N-[l\!-\!1]_qq^{N\!-\!1}+ 
\frac{(q^{- 2}\!-\!1)[l\!-\!1]_q}{q^{- 
1}\!+\!q^{N\!+\!1\!-\!2l}}  \right]\nn 
&=&\frac{[N\!-\!l\!+\!1]_q}{[l]_q} 
\frac{q^{\frac N2\!-\!l}\!+\!q^{l\!-\!\frac N2}} 
{q^{\frac N2\!+\!1\!-\!l}\!+\!q^{l\!-\!1\!-\!\frac N2}}  \qquad \qquad 
\mbox{if }H=U_qso(N).\label{tool2} 
\eea 
By the definition (\ref{defqepsilon}) of the $\varepsilon$-tensor and
(\ref{xixirel}) the claim is manifestly  true for $l=N$,
$$
{\cal P}^{-,N}{}^{i_1...i_N}_{j_1...j_N}=d_N \varepsilon^{i_1\ldots  i_N}
 \varepsilon^{j_1\ldots  j_N},
$$
because the rhs fulfills all
conditions (\ref{Plproj1}-\ref{Plproj3}). We prove the claim for
the remaining $l<N$ by induction, with $N$ inductive steps. Assume the
claim is true for $l=m\!+\!1$:
$$
{\cal P}^{-,{m\!+\!1}}{}^{i_1...i_{m\!+\!1}}_{j_1...j_{m\!+\!1}}=d_{m\!+\!1}  
U^{k_1}_{j_1}...U^{k_{m\!+\!1}}_{j_{m\!+\!1}} \varepsilon^{i_{m\!+\!2}\ldots
i_Ni_1\ldots i_{m\!+\!1}}  \varepsilon^{i_{m\!+\!2}\ldots i_Nk_1\ldots
k_{m\!+\!1}}. 
$$
Multiplying both sides by $U_{i_{m\!+\!1}}^{j_{m\!+\!1}}$ (and summing
of course also on the repeated indices $i_{m\!+\!1},j_{m\!+\!1}$) we find
on one hand
\bea
&&\left[\mbox{tr}_{m\!+\!1}\Big(U_{m\!+\!1}{\cal
P}^{-,{m\!+\!1}}\Big)\right]^{i_1...i_m}_{j_1...j_m} \nn
&&\qquad \qquad \qquad=d_{m\!+\!1}  
U^{k_1}_{j_1}...U^{k_m}_{j_m}U^2{}^{k_{m\!+\!1}}_{i_{m\!+\!1}}
\varepsilon^{i_{m\!+\!2}\ldots i_Ni_1\ldots i_{m\!+\!1}} 
\varepsilon^{i_{m\!+\!2}\ldots i_Nk_1\ldots k_{m\!+\!1}} \nn
&&\qquad\qquad \qquad\stackrel{(\ref{cycl_eps})}{=} d_{m\!+\!1}  
U^{k_1}_{j_1}...U^{k_m}_{j_m}
\varepsilon^{i_{m\!+\!1}\ldots i_Ni_1\ldots i_m} 
\varepsilon^{i_{m\!+\!}\ldots i_Nk_1\ldots k_m} \nonumber,
\eea
and on the other
\bea
\mbox{tr}_{m\!+\!1}\Big(U_{m\!+\!1}{\cal
P}^{-,{m\!+\!1}}\Big) &\stackrel{(\ref{ansatz1})}{=}& 
\mbox{tr}_{m\!+\!1}\Big(U_{m\!+\!1}
{\cal P}^{-,m}_{12...m}M^{-,m\!+\!1}_{m(m\!+\!1)}{\cal
P}^{-,m}_{12...l}\Big) \nn
&\stackrel{(\ref{tool1})}{=}& b_{m\!+\!1}
{\cal P}^{-,m}_{12...m},\nonumber
\eea
whence by comparison the claim for $l=m$ follows, 
because $d_mb_{m\!+\!1}=d_{m\!+\!1}$.

\subsubsection*{Proof of Proposition \ref{Hodge2} } 
\bea  
&& lhs(\ref{defHodge2}) \stackrel{(\ref{defthetaa'})}{=}   
{}^*[ \varphi(q^{\eta}S {\cal L}^{-,}{}^{i_1}_{j_1})\theta^{j_1}  
... \varphi(q^{\eta}S {\cal L}^{-,}{}^{i_p}_{j_p})\theta^{j_p}]\nn &&  
\stackrel{(\ref{framecond'})}{=}  
{}^*[ \varphi(q^{\eta}S {\cal L}^{-,}{}^{i_1}_{j_1})  ... 
\varphi(q^{\eta}S {\cal L}^{-,}{}^{i_p}_{j_p})\theta^{j_1}  
...\theta^{j_p}]\nn &&  
\stackrel{(\ref{Hbil})}{=}  
\varphi\Big(q^{p\eta}S({\cal L}^-{}^{i_p}_{j_p} ... 
 {\cal L}^{-,}{}^{i_1}_{j_1})\Big) {}^*[\theta^{j_1}  
...\theta^{j_p}]\nn &&  
\stackrel{(\ref{defHodge1})}{=} c_p q^{-Np/2}\Lambda^p  
\varphi\Big(S({\cal L}^-{}^{i_p}_{j_p} ... 
 {\cal L}^{-,}{}^{i_1}_{j_1})\Big) \theta^{h_{p+1}}...\theta^{h_N}  
\varepsilon_{h_N...h_{p+1}}{}^{j_1...j_p}\nn &&  
= c_p q^{-Np/2}\Lambda^p  
\varphi\Big( {\cal L}^{-,}{}^{l_{p+1}}_{j_{p+1}} ... 
{\cal L}^{-,}{}^{l_N}_{k_N} \Big)  
\varphi\Big(S({\cal L}^-{}^{i_p}_{j_p} ... 
 {\cal L}^{-,}{}^{i_1}_{j_1} {\cal L}^{-,}{}^{k_{p+1}}_{j_{p+1}}  
 ... {\cal L}^{-,}{}^{k_N}_{j_N}\Big)\nn &&  
\qquad\qquad \theta^{h_{p+1}}...\theta^{h_N}  
g_{h_Nl_N}...g_{ h_{p+1} l_{p+1}}  
\varepsilon^{j_N...j_{p+1}j_1...j_p}\nn &&  
\stackrel{(\ref{L-det})}{=} \!c_p q^{-Np/2}\Lambda^p  
\varphi\!\Big( \!{\cal L}^{-,}{}^{l_{p+1}}_{j_{p+1}} ... 
 {\cal L}^{-,}{}^{l_N}_{k_N}\! \Big)  
\theta^{h_{p+1}}...\theta^{h_N}  
\nn &&   
\qquad \qquad g_{h_Nl_N}...g_{ h_{p+1} l_{p+1}}  
\varepsilon^{k_N...k_{p+1}i_1...i_p}\nn &&  
\stackrel{(\ref{antipode})}{=} c_p q^{-Np/2}\Lambda^p  
g_{l_Nj_N}...g_{j_{p+1}l_{p+1}}  
\varphi\Big((S {\cal L}^{-,}{}^{l_{p+1}}_{h_{p+1}}) ... 
(S {\cal L}^{-,}{}^{l_N}_{h_N}) \Big)  
\nn &&  
\qquad \qquad\theta^{h_{p+1}}...\theta^{h_N}  
\varepsilon^{k_N...k_{p+1}i_1...i_p}\nn &&  
\stackrel{(\ref{defthetaa'})}{=} c_p q^{-N(p-N/2)}\Lambda^{2p-N}  
g_{l_Nk_N}...g_{l_{p+1} k_{p+1}}  
\xi^{l_{p+1}}...\xi^{l_N}  
\varepsilon^{k_N...k_{p+1}i_1...i_p}\nn &&  
= c_p q^{-N(p-N/2)}\Lambda^{2p-N}  
\xi^{l_{p+1}}...\xi^{l_N}  
\varepsilon_{l_N...l_{p+1}}{}^{i_1...i_p}\nn &&  
= rhs(\ref{defHodge2})  \nonumber
\eea  
  
\subsection*{Proof of Proposition \ref{lapl}}

We now evaluate  the lhs(\ref{Laplacian}) on each 
${\cal DC}^p$. We  find
\bea
&&(d\,{}^*\, d\,{}^*+{}^*\, d\,{}^*d)\,\1={}^*\, d\,{}^*d\,\1=
{}^*\, d\,{}^*\xi^{i_1}\partial_{i_1}
\nn && =q^{-N(1-\frac N2)}c_1{}^*\,
d\,\xi^{i_2}...\xi^{i_N}\varepsilon_{i_N...i_2}{}^{i_1}\Lambda^{2\!-\!N}\partial_{i_1}
\nn && = (-1)^{N-1}q^{-N(1-\frac N2)}c_1{}^*\,
\xi^{i_2}...\xi^{i_N}d\varepsilon_{i_N...i_2}{}^{i_1}\Lambda^{2\!-\!N}\partial_{i_1}
\nn && =(-1)^{N-1}q^{-N}c_1c_N\varepsilon^{i_2...i_Nj}\Lambda^N\partial_j
\varepsilon_{i_N...i_2}{}^{i_1}\Lambda^{2\!-\!N}\partial_{i_1}
\nn && = q^2c_1c_N U^{-1}{}^j_h\varepsilon^{hi_2...i_N}\partial_j
\varepsilon_{i_N...i_2}{}^{i_1}\partial_{i_1}\Lambda^2
\nn && =q^2c_1c_N g^{lj}\varepsilon_l{}^{i_2...i_N}\partial_j
\varepsilon_{i_N...i_2}{}^{i_1}\partial_{i_1}\Lambda^2
\nn && =\frac{q^2c_1c_N}{d_1} g^{i_1j}\partial_j\partial_{i_1}\Lambda^2=
\frac{q^2c_1c_N}{d_1} \Box\Lambda^2
\nn && =\frac{q^2c_1}{c_0[N]_q}\frac{q^{-\frac N2}\!+\!q^{\frac N2}}
{q^{1\!-\!\frac N2}\!+\!q^{\frac N2\!-\!1}}\,\Box\Lambda^2\nonumber
\eea
for $p=0$, for $p=N$
\bea
&&(d\,{}^*\, d\,{}^*+{}^*\, d\,{}^*d)\,d^N\!x=
d\,{}^*\, d\,{}^*\,d^N\!x=
q^{-\frac {N^2}2}c_Nd\,{}^*\, d\Lambda^N 
=q^{-\frac{N^2}2}c_Nd\,{}^*\,\xi^{i_1}\partial_{i_1}\Lambda^N 
\nn && =q^{-N}c_Nc_1d\,\xi^{i_2}...\xi^{i_N}\,
\varepsilon_{i_N...i_2}{}^{i_1}\Lambda^{2\!-\!N} \partial_{i_1}\Lambda^N 
\nn && =(-1)^{N\!-\!1}q^{2\!-\!2N}c_Nc_1\xi^{i_2}...\xi^{i_N}d\,
\varepsilon_{i_N...i_2}{}^{i_1}\,\partial_{i_1}\Lambda^2
\nn && =(-1)^{N\!-\!1}q^{2\!-\!2N}c_Nc_1\xi^{i_2}...\xi^{i_N}\xi^j\partial_j
\varepsilon_{i_N...i_2}{}^{i_1}\,\partial_{i_1}\Lambda^2
\nn && =(-1)^{N\!-\!1}q^{2\!-\!2N}c_Nc_1\varepsilon^{i_2...i_Nj}\,d^N\!x\,
\varepsilon_{i_N...i_2}{}^{i_1}\,\partial_j\partial_{i_1}\Lambda^2
\nn && =q^2c_Nc_1g^{lj}\varepsilon_l{}^{i_2...i_N}
\varepsilon_{i_N...i_2}{}^{i_1}\,\partial_j\partial_{i_1}\Lambda^2d^N\!x
\nn && =\frac{q^2c_1c_N}{d_1} g^{i_1j}\partial_j\partial_{i_1}\Lambda^2\,d^N\!x=
\frac{q^2c_1c_N}{d_1} \Box\Lambda^2d^N\!x
\nn && =\frac{q^2c_1}{c_0[N]_q}\frac{q^{-\frac N2}\!+\!q^{\frac N2}}
{q^{1\!-\!\frac N2}\!+\!q^{\frac N2\!-\!1}}\,\Box\Lambda^2\,d^N\!x\nonumber
\eea
for $p=N$, whereas for $p=1,2,...,N-1$ we find on one hand
\bea 
&&d\,{}^*\, d\,{}^*\xi^{i_1}...\xi^{i_p} 
\stackrel{(\ref{defHodge2})}{=}(-1)^{N\!-\!p} 
q^{-N(p\!-\!\frac N2)+N-2p} 
c_p\, d\,{}^*\, \xi^{i_{p+1}}...\xi^{i_N} 
\varepsilon_{i_N...i_{p+1}}{}^{i_1...i_p}\Lambda^{2p\!-\!N}\,d \nn 
&&\stackrel{(\ref{defHodge2})}{=}(-1)^{N\!-\!p} q^{-N(p\!-\!\frac N2) 
+N-2p-N(N\!-\!p\!+\!1\!-\!\frac N2)} c_p c_{N\!-\!p\!+\!1} 
\, d\, \nn 
&&\qquad \qquad\xi^{h_1}...\xi^{h_{p-1}} 
\varepsilon_{h_{p-1}...h_1}{}^{i_{p+1}...i_Ni} 
\varepsilon_{i_N...i_{p+1}}{}^{i_1...i_p}\Lambda^2\partial_i \nn 
&&\stackrel{(\ref{cycl_eps} )}{=}(-1)^{N\!-\!1} q^{-2p}c_p 
c_{N\!-\!p\!+\!1}  \xi^{h_1}...\xi^{h_{p-1}}d 
(-1)^{N\!-\!1}\nn 
&&\qquad g^{l_pi}\varepsilon_{l_ph_{p-1}...h_1}{}^{i_{p+1}...i_N} 
\varepsilon_{i_N...i_{p+1}}{}^{i_1...i_p}\Lambda^2\partial_i \nn 
&&\stackrel{(\ref{Pa_eps_rel})}{=}q^{2-2p} 
\frac{c_p c_{N\!-\!p\!+\!1}}{d_p}\,
 \xi^{h_1}...\xi^{h_{p\!-\!1}}d\, 
{\cal P}^{a,p}{}^{i_1...i_{p-1}i_p}_{h_1...h_{p\!-\!1}j_p}  
\partial^{j_p} \Lambda^2\nonumber 
\eea
and on the other
\bea 
&&{}^*\, d\,{}^*\, d\xi^{i_1}...\xi^{i_p}= 
(-1)^p\,{}^*\, d\,{}^*\xi^{i_1}...\xi^{i_p}d\nn 
&&\stackrel{(\ref{defHodge2})}{=}(-1)^pc_{p\!+\!1} 
q^{-N(p\!+\!1\!-\!\frac N2)} 
\,{}^*\,d\,\xi^{i_{p\!+\!2}}...\xi^{i_N} 
\varepsilon_{i_N...i_{p+2}}{}^{i_1...i_{p+1}} 
\Lambda^{2p\!+\!2\!-\!N}\partial_{i_{p+1}}  
\nn &&=(-1)^{N\!-\!1}c_{p\!+\!1}q^{-N(p\!-\!\frac N2)\!-\!2p\!-\!2} 
\,{}^*\,\xi^{i_{p\!+\!2}}...\xi^{i_N}\xi^{j_N} 
\varepsilon_{i_N...i_{p+2}}{}^{i_1...i_{p+1}} 
\Lambda^{2p\!+\!2\!-\!N}\partial_{j_N}\partial_{i_{p+1}}  
\nn &&\stackrel{(\ref{defHodge2})}{=} 
(-1)^{N\!-\!1}c_{p\!+\!1}c_{N\!-\!p} 
q^{-\!2p\!-\!2} \,\xi^{h_1}...\xi^{h_p} 
\varepsilon_{h_p...h_1}{}^{i_{p\!+\!2}...i_Nj_N} 
\varepsilon_{i_N...i_{p+2}}{}^{i_1...i_{p+1}} 
\Lambda^2\partial_{j_N}\partial_{i_{p+1}}  
\nn &&\stackrel{(\ref{cycl_eps})}{=} 
c_{p\!+\!1}c_{N\!-\!p} 
q^{-\!2p\!-\!2} \,\xi^{h_1}...\xi^{h_p} g^{h_{p\!+\!1}j_N} 
\varepsilon_{h_{p\!+\!1}h_p...h_1}{}^{i_{p\!+\!2}...i_N} 
\varepsilon_{i_N...i_{p+2}}{}^{i_1...i_{p+1}} 
\Lambda^2\partial_{j_N}\partial_{i_{p+1}}  
\nn &&\stackrel{(\ref{Pa_eps_rel})}{=} 
q^{2\!-\!2p}\frac{c_{p\!+\!1}c_{N\!-\!p}}{d_{p\!+\!1}}
\,\xi^{h_1}...\xi^{h_p}  {\cal P}^{a,p\!+\!1}
{}^{i_1...i_{p\!+\!1}}_{h_1...h_{p\!+\!1}}    g^{h_{p\!+\!1}j_N}
 \partial_{j_N}\partial_{i_{p+1}} \Lambda^2   \nn
&&\stackrel{(\ref{ansatz1})}{=}  \frac{q^{2\!-\!2p}}{[p\!+\!1]_q}
\frac{c_{p\!+\!1}c_{N\!-\!p}}{d_{p\!+\!1}}
\,\xi^{h_1}...\xi^{h_p}  {\cal
P}^{a,p}{}^{i_1...i_{p\!-\!1}i_p}_{h_1...h_{p\!-\!1}j_p} \left[ q^p
\delta^{j_p}_{h_p}g^{i_{p\!+\!1}j_N}\right.\nn
&&\qquad\quad \left .\! -\! [p]_q \hat R^{j_pi_{p\!+\!1}}_{h_ph_{p\!+\!1}}  
 g^{h_{p\!+\!1}j_N}\!-\!  \frac{k[p]_q}{1\!+\!q^{N-2p}}
\delta^{j_N}_{h_p}g^{j_pi_{p\!+\!1}}\right] \partial_{j_N}\partial_{i_{p+1}}
\Lambda^2  
\nn &&\stackrel{(\ref{gRrel})}{=}  \frac{q^{2\!-\!2p}}{[p\!+\!1]_q}
\frac{c_{p\!+\!1}c_{N\!-\!p}}{d_{p\!+\!1}}
\,\xi^{h_1}...\xi^{h_p}  {\cal
P}^{a,p}{}^{i_1...i_{p\!-\!1}i_p}_{h_1...h_{p\!-\!1}j_p} \left[ q^p
\delta^{j_p}_{h_p}g^{i_{p\!+\!1}j_N}\right.\nn
&&\qquad\quad \left .\! -\! [p]_q
 g^{j_ph_{p\!+\!1}}\hat R^{-1}{}^{i_{p\!+\!1}j_N}_{h_{p\!+\!1}h_p}   \!-\! 
\frac{k[p]_q}{1\!+\!q^{N-2p}}\delta^{j_N}_{h_p}g^{j_pi_{p\!+\!1}}\right]
\partial_{j_N}\partial_{i_{p+1}} \Lambda^2   \nn
&&\stackrel{(\ref{ddrel}),(\ref{projectorR})}{=} 
\frac{q^{2\!-\!2p}}{[p\!+\!1]_q} \frac{c_{p\!+\!1}c_{N\!-\!p}}{d_{p\!+\!1}}
\,\xi^{h_1}...\xi^{h_p}  {\cal
P}^{a,p}{}^{i_1...i_{p\!-\!1}i_p}_{h_1...h_{p\!-\!1}j_p} \left[
\delta^{j_p}_{h_p}g^{i_{p\!+\!1}j_N}\Big(q^p\! -\! [p]_q\frac
k{\mu}\Big)\right. \nn &&\qquad\quad \left . \!-\!  [p]_q
\left(q^{-1}+\frac
k{1\!+\!q^{N-2p}}\right) \delta^{j_N}_{h_p}g^{j_pi_{p\!+\!1}}\right]
\partial_{j_N}\partial_{i_{p+1}} \Lambda^2    \nn
&&= \frac{q^{2\!-\!2p}}{[p\!+\!1]_q} \frac{c_{p\!+\!1}c_{N\!-\!p}}{d_{p\!+\!1}}
\left[\frac{q^{p\!+\!1\!-\!\frac N2}\!+\!q^{\frac N2\!-\!p\!-\!1}}
{q^{1\!-\!\frac N2}\!+\!q^{\frac N2\!-\!1}}\,\xi^{i_1}...\xi^{i_p} \Box
\right. \nn &&\qquad\quad \left . \!-\!  [p]_q
\frac{q^{p\!+\!1\!-\!\frac N2}\!+\!q^{\frac N2\!-\!p\!-\!1}}
{q^{p\!-\!\frac N2}\!+\!q^{\frac N2\!-\!p}}\xi^{h_1}...\xi^{h_{p\!-\!1}}d \,
{\cal P}^{a,p}{}^{i_1...i_{p\!-\!1}i_p}_{h_1...h_{p\!-\!1}j_p}
\partial^{j_p}\right]\Lambda^2                    \label{interm}
\eea 
In order that the second term in the square bracket be
opposite of $d\,{}^*\, d\,{}^*\xi^{i_1}...\xi^{i_p} $ it must be
\bea
0&\stackrel{!}{=}&\frac{c_p c_{N\!-\!p\!+\!1}}{d_p}-
\frac{q^{p\!+\!1\!-\!\frac N2}\!+\!q^{\frac N2\!-\!p\!-\!1}}
{q^{p\!-\!\frac N2}\!+\!q^{\frac N2\!-\!p}}
\frac{[p]_q}{[p\!+\!1]_q}\frac{c_{p\!+\!1}c_{N\!-\!p}}{d_{p\!+\!1}}
\nn &=&\frac{q^{-\frac N2}\!+\!q^{\frac N2}}
{q^{p\!-\!\frac N2}\!+\!q^{\frac N2\!-\!p}}
\frac{[p]_q![N\!-\!p\!-\!1]_q!}{[N]_q!d_N}\left\{[N\!-\!p]_q c_p
c_{N\!-\!p\!+\!1} - [p]_q c_{p\!+\!1}c_{N\!-\!p}\right\}
\nonumber
\eea
namely, for $p=1,2,...,N\!-\!1$ 
$$
[N\!-\!p]_q c_p
c_{N\!-\!p\!+\!1} - [p]_q c_{p\!+\!1}c_{N\!-\!p}=0.
$$
This recursion relation is solved by
\be
c_{p\!+\!1}c_{N\!-\!p}
={[N\!-\!1]_q\choose{[p]_q}}\,c_1c_N.     \label{solved}
\ee
When replaced in (\ref{interm}) we find, on all
of ${\cal DC}^*$, and in particular on $\Omega^*$,
\be
(d\,{}^*\, d\,{}^*+{}^*\, d\,{}^*d)\,\xi^{i_1}...\xi^{i_p}
=\frac{q^2c_1}{c_0[N]_q}\frac{q^{-\frac N2}\!+\!q^{\frac N2}}
{q^{1\!-\!\frac N2}\!+\!q^{\frac N2\!-\!1}}\,\Box\Lambda^2
\xi^{i_1}...\xi^{i_p}.                  \label{Lap}
\ee
Taking $\prod\limits_{p=1}^{N-2}$ of both sides 
of (\ref{solved}) and multiplying the
result by $c_0(c_1c_N)^2$  we obtain
$$
d_1d_2...d_Nc_N=c_0(c_1c_N)^N{[N\!-\!1]_q\choose{[1]_q}}...
{[N\!-\!1]_q\choose{[N\!-\!2]_q}},
$$
implying, because of (\ref{defd_p}),
$$
\left(\frac{c_1}{c_0[N]_q}\right)^N=
\left(\frac{c_1c_N}{d_N[N]_q)}\right)^N=
\frac{c_N}{c_0\,[N]_q!}\prod\limits_{l=0}^{N-1}
\frac{q^{l\!-\!\frac N2}\!+\!q^{\frac N2\!-\!l}} 
{q^{-\frac N2}\!+\!q^{\frac N2}}.
$$
Before proceeding we note that we are still free to choose the value
of $d_N=c_0c_N$ and the normalization of $c_N$ w.r.t. $c_0$,
in other words we are free to choose the values of both $c_0,c_N$.
We choose
\be
c_N=1,\qquad\qquad\qquad c_0=d_N=
\frac 1{[N]_q!}\prod\limits_{l=0}^{N-1}
\frac{q^{l\!-\!\frac N2}\!+\!q^{\frac N2\!-\!l}} 
{q^{1\!-\!\frac N2}\!+\!q^{\frac N2\!-\!1}};       \label{choice}
\ee
the first choice guarantees that ${}^*\1=dV$, and therefore
also ${}^*dV=\1$ in view of $*^2$. As a consequence,
$$
\left(\frac{c_1}{c_0[N]_q}\right)^N=
\left(\frac{q^{1\!-\!\frac N2}\!+\!q^{\frac N2\!-\!1}} 
{q^{-\frac N2}\!+\!q^{\frac N2}}\right)^N
$$
implying
$$
c_1=c_0[N]_q\frac{q^{1\!-\!\frac N2}\!+\!q^{\frac N2\!-\!1}} 
{q^{-\frac N2}\!+\!q^{\frac N2}}=
\frac 1{[N\!-\!1]_q!}\prod\limits_{l=1}^{N-1}
\frac{q^{l\!-\!\frac N2}\!+\!q^{\frac N2\!-\!l}} 
{q^{1\!-\!\frac N2}\!+\!q^{\frac N2\!-\!1}}.
$$
Multiplying both sides of (\ref{solved}) by 
$1/c_{N\!-\!p}=c_p/d_p$ and using (\ref{defd_p}) we 
thus obtain the
recursion relation
$$
c_{p\!+\!1}=c_p [N\!-\!p]_q
\frac{q^{1\!-\!\frac N2}\!+\!q^{\frac N2\!-\!1}}
{q^{p\!-\!\frac N2}\!+\!q^{\frac N2\!-\!p}}.
$$
Solving the latter, the final result for $p=0,1,...,N$ reads (\ref{c_p=}).
As for (\ref{Lap}), we find that on all
of ${\cal DC}^*$, and in particular on $\Omega^*$, (\ref{Laplacian}) holds.

Finally, we see that the normalization condition (\ref{choice})
for $d_N=d_0$ implies a specific value for the normalization
constant 
$$
\gamma_N\equiv\cases{ (\varepsilon^{12...N})^{-1}\cr
(\varepsilon^{-n(1\!-\!n)...n})^{-1}}
$$ 
in (\ref{defdNx}), which can be computed case by case.
In particular, for the cases $N=3,4$ this implies a
different normalization w.r.t. the one adopted in section \ref{eps}.

\end{document}